\documentclass[format=acmsmall,screen,timestamp=false,a4paper]{acmart}
\pdfoutput=1
\newif\ifarxiv
\arxivtrue

\ifarxiv
\setkeys{acmart.cls}{nonacm}
\else
\setkeys{acmart.cls}{draft}
\fi
\setcopyright{none}
\acmDOI{}

\acmJournal{TOMS}

\usepackage{amsmath}
\usepackage{caption}
\usepackage{subcaption}
\usepackage{pgfplots}
\usepackage{pgfplotstable}
\pgfplotscreateplotcyclelist{mycolor}{
  red,every mark/.append style={fill=red},mark=*\\
  blue,every mark/.append style={fill=blue},mark=*\\
  orange,every mark/.append style={fill=orange},mark=*\\
  red,densely dotted, every mark/.append style={fill=red},mark=square*\\
  blue,densely dotted, every mark/.append style={fill=blue},mark=square*\\
  orange,densely dotted every mark/.append style={fill=orange},mark=square*\\
}

\usetikzlibrary{3d,calc}
\tikzset{persp/.style={scale=0.5,x={(-0.8cm,-0.4cm)},y={(0.8cm,-0.4cm)}, z={(0cm,1cm)}},points/.style={fill=white,draw=black,thick}}

\usepackage{siunitx}
\usepackage{listings}
\usepackage{xspace}
\usepackage{booktabs}
\usepackage{multirow}
\DeclareMathOperator{\Div}{div}
\DeclareMathOperator{\curl}{curl}

\newcommand{\R}{\mathbb{R}}

\newcommand{\calP}{\mathcal{P}}
\newcommand{\calQ}{\mathcal{Q}}
\newcommand{\calS}{\mathcal{S}}
\newcommand{\hcurl}{\ensuremath{{H}(\curl)}\xspace}
\newcommand{\hdiv}{\ensuremath{{H}(\Div)}\xspace}

\newcommand{\cancel}[1]{}

\title[Trimmed Serendipity elements in Firedrake]{Bringing Trimmed Serendipity Methods to Computational Practice in Firedrake}

\author{Justin Crum}
\email{jcrum@math.arizona.edu}
\affiliation{%
  \institution{University of Arizona}
  \department{Department of Mathematics}
  \city{Tucson}
  \state{AZ}
  \country{USA}
}

\author{Cyrus Cheng}
\email{cyrus.cheng15@alumni.imperial.ac.uk}
\affiliation{%
  \institution{Imperial College London}
  \department{Department of Mathematics}
  \city{London}
  \postcode{SW7 2AZ}
  \country{UK}}

\author{David A. Ham}
\affiliation{%
  \institution{Imperial College London}
  \department{Department of Mathematics}
  \city{London}
  \postcode{SW7 2AZ}
  \country{UK}}
\email{david.ham@imperial.ac.uk}

\author{Lawrence Mitchell}
\orcid{0000-0001-8062-1453}
\affiliation{%
  \institution{Durham University}
  \department{Department of Computer Science}
  \streetaddress{Upper Mountjoy}
  \city{Durham}
  \postcode{DH1 3LE}
  \country{UK}}
\email{lawrence.mitchell@durham.ac.uk}

\author{Robert C. Kirby}
\affiliation{%
  \institution{Baylor University}
  \department{Department of Mathematics}
  \streetaddress{1410 S.~4th Street}
  \city{Waco}
  \state{TX}
  \country{USA}
}
\email{robert\_kirby@baylor.edu}

\author{Joshua A. Levine}
\affiliation{%
  \institution{University of Arizona}
  \department{Department of Computer Science}
  \city{Tucson}
  \state{AZ}
  \country{USA}
}
\email{josh@email.arizona.edu}

\author{Andrew Gillette}
\affiliation{%
  \institution{University of Arizona}
  \department{Department of Mathematics}
  \city{Tucson}
  \state{AZ}
  \country{USA}}
\email{agillette@math.arizona.edu}

\ccsdesc[300]{Mathematics of computing~Mathematical software performance}
\ccsdesc[100]{Mathematics of computing~Solvers}
\ccsdesc[500]{Mathematics of computing~Partial differential equations}
\ccsdesc[500]{Mathematics of computing~Discretization}

\begin{abstract}
       We present an implementation of the trimmed serendipity finite element family, using the open source finite element package Firedrake.  The new elements can be used seamlessly within the software suite for problems requiring $H^1$, \hcurl, or \hdiv-conforming elements on meshes of squares or cubes.  To test how well trimmed serendipity elements perform in comparison to traditional tensor product elements, we perform a sequence of numerical experiments including the primal Poisson, mixed Poisson, and Maxwell cavity eigenvalue problems.  Overall, we find that the trimmed serendipity elements converge, as expected, at the same rate as the respective tensor product elements while being able to offer significant savings in the time or memory required to solve certain problems.
\end{abstract}

\usepackage[capitalise]{cleveref}
\begin{document}
  \maketitle

  \section{Introduction}

    Over the last 15 years, conforming finite element methods for Hodge-Laplace type problems on simplicial and cubical meshes have been analyzed and categorized using the mathematical theory of Finite Element Exterior Calculus (FEEC)~\cite{AFW2006,AFW2010,ABB2012}. 
    While this effort was initially focused on placing related theoretical results under a single, unified mathematical framework, it also exposed the frontier of knowledge and spawned a wealth of spin-off projects that improved awareness, understanding, and implementation of the methods described.
    Among these follow-on projects was a revisiting of the notion of ``serendipity'' finite elements, a particular type of finite element method dating back to the 1970s.
    Named for its seemingly too-good-to-be-true computational benefit, a serendipity finite element method converged to the correct solution of a partial differential equation (PDE) at an equivalent rate but with fewer degrees of freedom than the corresponding tensor product method.
    The simplest and most well-known serendipity method replaces the quadratic 9-node square element with a ``serendipity'' 8-node element that has no interior degree of freedom, but still converges at a quadratic rate in the appropriate sense.
    
    The core idea of serendipity elements, i.e.\ reducing degrees of freedom from tensor product methods without reducing the accuracy of the method, floated around computational engineering communities until it found a resurgence under the FEEC framework.
    In a seminal work by Arnold and Awanou~\cite{arnold2011serendipity} the scalar-valued serendipity elements were formalized as providing approximation on a space of polynomials $\mathcal{S}_r$ that nests between maximum-degree $r$ polynomials $\mathcal{P}_r$ and tensor product degree $r$ polynomials $\mathcal{Q}_r$.  
    In a subsequent work~\cite{arnold2014finite}, this notion was extended using FEEC to introduce ``serendipity vector elements'' to the literature -- analogues of the famous N\'ed\'elec elements~\cite{N1980,N1986} -- which led to the widely circulated Periodic Table of Finite Elements \protect{\cite{arnold2014periodic}}.
    
    Despite excitement at these developments, implementing serendipity elements and realizing any potential computational benefits proved to be a significant challenge.
    While creating data structures for ``all polynomials up to degree $r$'' is straightforward, simply trying to write down the approximation spaces used by serendipity elements -- especially the vector-valued elements -- requires a substantial amount of mathematical notation and explanation.
    As a consequence, the serendipity elements have never been implemented or rigorously studied beyond a few special use cases.
    
    Recent work on a closely related family of methods -- called \textit{trimmed serendipity elements}~\cite{gillette2019trimmed} -- has opened the door to potential widespread usage and benefit from the notions of serendipity theory just described.
    The trimmed serendipity spaces are identical to the serendipity spaces in the scalar-valued cases but are distinct in the vector-valued cases where, notably, they have fewer degrees of freedom than the serendipity elements of the same order.
    Since the computational motivation for serendipity methods is entirely about reducing degree of freedom count while preserving approximation power, the smaller dimensionality of the trimmed serendipity spaces obviates the need to implement the vector-valued ``non-trimmed'' serendipity elements.
    
    Implementation of trimmed serendipity spaces is feasible by merging two independent efforts: the systematic definition of computational basis functions for these methods from work by Gillette, Kloefkorn and Sanders~\cite{gillette2019computational} and the easily extensible open-source finite element software package Firedrake ~\cite{rathgeber2016firedrake}.
    The Unified Form Language \cite{Logg:2012,alnaes2014unified}, inspired by FEEC, provides a common backbone for translating the basis functions from their formal statement into a code structure using a well-established, high level interface.
    We explain the idea of the implementation by walking through the process of discretizing a PDE and selecting a finite element method for approximating its solution.

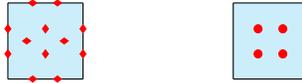
\begin{figure}[htbp]
\begin{tikzpicture}

\filldraw [cyan, opacity=0.2] (0,0,1) -- (1,0,1) -- (1,1,1) -- (0,1,1) -- cycle;

\draw (0,1,1) -- (1,1,1);

\draw (0,0,1) -- (0,1,1);
\draw (0,0,1) -- (1,0,1);
\draw (1,0,1) -- (1,1,1);

\node[mark size=1.5pt,color=red, rotate=0] at (0.5,0.33,1) {\pgfuseplotmark{diamond*}};
\node[mark size=1.5pt,color=red, rotate=0] at (0.5,0.66,1) {\pgfuseplotmark{diamond*}};
\node[mark size=1.5pt,color=red, rotate=90] at (0.25,0.5,1) {\pgfuseplotmark{diamond*}};
\node[mark size=1.5pt,color=red, rotate=90] at (0.75,0.5,1) {\pgfuseplotmark{diamond*}};


\node[mark size=1.5pt,color=red, rotate=270] at (0.33,1,1) {\pgfuseplotmark{diamond*}};
\node[mark size=1.5pt,color=red, rotate=90] at (0.66,1,1) {\pgfuseplotmark{diamond*}};
\node[mark size=1.5pt,color=red, rotate=90] at (0.66,0,1) {\pgfuseplotmark{diamond*}};
\node[mark size=1.5pt,color=red, rotate=90] at (0.33,0,1) {\pgfuseplotmark{diamond*}};
\node[mark size=1.5pt,color=red, rotate=0] at (0,0.66,1) {\pgfuseplotmark{diamond*}};
\node[mark size=1.5pt,color=red, rotate=0] at (0,0.33,1) {\pgfuseplotmark{diamond*}};
\node[mark size=1.5pt,color=red, rotate=0] at (1,0.66,1) {\pgfuseplotmark{diamond*}};
\node[mark size=1.5pt,color=red, rotate=0] at (1,0.33,1) {\pgfuseplotmark{diamond*}};
\filldraw [cyan, opacity=0.2] (3,0,1) -- (4,0,1) -- (4,1,1) -- (3,1,1) -- cycle;

\draw (3,1,1) -- (4,1,1);
\draw (3,0,1) -- (3,1,1);
\draw (3,0,1) -- (4,0,1);
\draw (4,0,1) -- (4,1,1);

\filldraw [red] (3.33,0.66,1) circle (1.5pt);
\filldraw [red] (3.33,0.33,1) circle (1.5pt);
\filldraw [red] (3.66,0.33,1) circle (1.5pt);
\filldraw [red] (3.66,0.66,1) circle (1.5pt);

\end{tikzpicture}
\caption{The degree 2 \texttt{RTCF} tensor product element (Raviart-Thomas \hcurl elements on quadrilaterals, \cite{raviart1977mixed}) (left) used in \cref{lst:mixedPoissonCode}.  The \texttt{RTCF} element is an example of a tensor product element in 2D, which depending upon the orientation of the DOFs on the edges, can be used for either \hdiv or \hcurl problems.  The right displays a similar example for the \texttt{DQ} element, which is an $L^2$-conforming tensor product element in 2D, and is necessary to form the stable pair for the mixed Poisson problem.\label{fig:RTCF}}
\Description{An example of what a reference RTCF finite element looks like.  Dots represent DOFs on the face or edges.}
\end{figure}

\paragraph{Motivating example}
 \cref{lst:mixedPoissonCode} provides a snippet of code that a user could write to approximate a solution to the mixed Poisson equation on the domain $\Omega := [0, 1] \times [0,1]$ with boundary $\Gamma$.  The formal problem statement of the continuous weak form in this case is: find $\sigma \in \Sigma :=~$\hdiv and $u \in V := L^2$ such that:

\begin{equation}\label{eq:MixedPoisson}
\begin{tabular}{rllll}

$\displaystyle\int_\Omega (\sigma \cdot \tau + \nabla\cdot u\tau) \text{ d}x$ & $=~~ 0$ & $\forall~ \tau\in\Sigma$, \\

$\displaystyle\int_\Omega \nabla \cdot \sigma v \text{ d}x$ & $= ~~\displaystyle -\int_\Omega fv \text{ d}x$ & $\forall~ v\in V$.
\end{tabular}
\end{equation}
We assume homogeneous Dirichlet boundary conditions so that the Dirichlet data in $u$ vanishes on $\Gamma$.
Solving the discretized version of (\ref{eq:MixedPoisson}) requires choosing a suitable pair of finite element spaces to create a stable method.  On a mesh of squares, the typical stable pair of tensor product finite elements would be \texttt{RTCF} and \texttt{DQ} for \hdiv and $L^2$ respectively.  These elements are visualized in \cref{fig:RTCF} and are part of the tensor product family of elements. 
We demonstrate the tensor product pairing in \cref{lst:mixedPoissonCode}, where the order of the vector and scalar elements are offset by 1 in accordance with the theory for optimal convergence rates.

\begin{lstlisting}[float=htpb,caption={Basic Firedrake implementation of the mixed Poisson problem showcasing where to choose the elements that are used and how to create the equations in Firedrake's notation.}, label={lst:mixedPoissonCode}, numbers=left, firstnumber=1, xleftmargin=20pt,  xrightmargin=20pt]
polyDegree = 2
numberOfCells = 2**5
mesh = UnitSquareMesh(numberOfCells, numberOfCells, quadrilateral=True)
hDivSpace = FunctionSpace(mesh, "RTCF", polyDegree)
l2Space = FunctionSpace(mesh, "DQ", polyDegree - 1)
mixedSpace = hDivSpace * l2Space

sigma, u = TrialFunctions(mixedSpace)
tau, v = TestFunctions(mixedSpace)

x, y = SpatialCoordinate(mesh)
uex = sin(pi*x)*sin(pi*y)

f = -div(grad(uex))
a = (dot(sigma, tau) + div(tau)*u + div(sigma)*v)*dx
l = -f*v*dx
w = Function(mixedSpace)
solve(a == l, w)
\end{lstlisting}

An important strength of Firedrake is its modular structure for both users and developers.
For the user, swapping to trimmed serendipity elements to solve the mixed Poisson problem is now only a matter of modifying lines 4--5 in \cref{lst:mixedPoissonCode} to the appropriate identifiers, \texttt{SminusDiv} and \texttt{DPC} inside the \texttt{FunctionSpace} calls that define \texttt{hDivSpace} and \texttt{l2Space}.
For developers, implementing a new element type -- such as trimmed serendipity -- is simply a matter of defining a suitable computational basis and connecting it to the intermediate interfaces in the included libraries.  

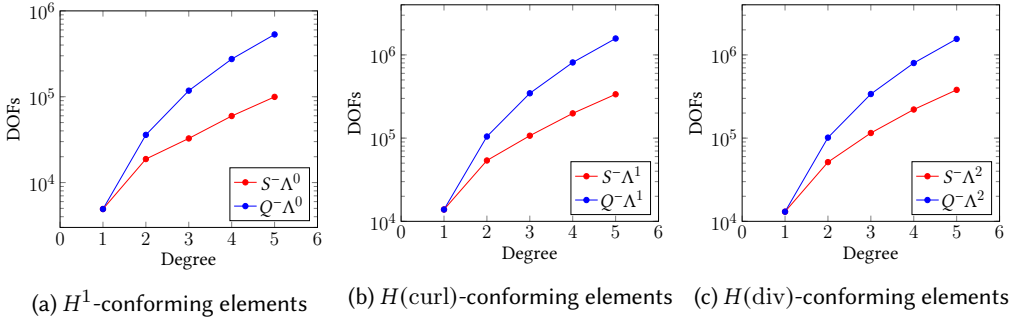
\begin{figure}[htbp]
  \centering
  \begin{subfigure}[h]{0.32\textwidth}
    \begin{tikzpicture}[scale=0.5]
      \begin{semilogyaxis}[
        xlabel={Degree}, ylabel={DOFs},
        ticklabel style={font=\LARGE},
        xlabel style={font=\LARGE},
        ylabel near ticks,
        ylabel style={font=\LARGE},
        ymax=1e6, ymin=3e3, xmax=6, xmin=0,
        legend pos=south east, legend style={font=\LARGE} ,
        cycle list name=mycolor]
        \addplot
        table [x=Degree,y=Dofs, col sep=comma]{csvs/SerendipityDegDofsTwo.csv};
        \addlegendentry{$S^- \Lambda^0$}
        \addplot
        table [x=Degree,y=Dofs, col sep=comma]{csvs/LagrangeDegDofs.csv};
        \addlegendentry{$Q^-\Lambda^0$}
      \end{semilogyaxis}
      \end{tikzpicture}
      \caption{$H^1$-conforming elements\label{fig:LagrangeDegDofs}}
    \end{subfigure}
    \begin{subfigure}[h]{0.32\textwidth}
      \begin{tikzpicture}[scale=0.5]
        \begin{semilogyaxis}[
          xlabel={Degree}, ylabel={DOFs},
          ticklabel style={font=\LARGE},
          xlabel style={font=\LARGE},
          ylabel near ticks,
          ylabel style={font=\LARGE},
          ymax=4e6, ymin=1e4, xmax=6, xmin=0,
          legend pos=south east, legend style={font=\LARGE} ,
          cycle list name=mycolor]
        \addplot
        table [x=Degree,y=Dofs, col sep=comma]{csvs/SminusCurlDegDofs.csv};
        \addlegendentry{$S^- \Lambda^1$}
        \addplot
        table [x=Degree,y=Dofs, col sep=comma]{csvs/NCEDegDofs.csv};
        \addlegendentry{$Q^-\Lambda^1$}
      \end{semilogyaxis}
      \end{tikzpicture}
      \caption{\hcurl-conforming elements\label{fig:CurlDegDofs}}
  \end{subfigure} 
  \begin{subfigure}[h]{0.32\textwidth}
    \begin{tikzpicture}[scale=0.5]
      \begin{semilogyaxis}[
        xlabel={Degree}, ylabel={DOFs},
        ticklabel style={font=\LARGE},
        xlabel style={font=\LARGE},
        ylabel near ticks,
        ylabel style={font=\LARGE},
        ymax=4e6, ymin=1e4, xmax=6, xmin=0,
        legend pos=south east, legend style={font=\LARGE} ,
        cycle list name=mycolor]
        \addplot
        table [x=Degree,y=Dofs, col sep=comma]{csvs/SminusDivDegDofs.csv};
        \addlegendentry{$S^- \Lambda^2$}
        \addplot
        table [x=Degree,y=Dofs, col sep=comma]{csvs/NCFDegDofs.csv};
        \addlegendentry{$Q^-\Lambda^2$}
      \end{semilogyaxis}
      \end{tikzpicture}
      \caption{\hdiv-conforming elements\label{fig:SminusDegDofs}}
    \end{subfigure}
  \caption{Comparison of the degrees of freedom (DOFs) required for a trimmed serendipity element $\mathcal{S}^- \Lambda^k$ vs a tensor product element $\mathcal{Q}^- \Lambda^k$, as calculated on a $[0,1]^3$ mesh with a total of $16^3$ cubes.  After order 1, where the element types coincide, the trimmed serendipity elements have strictly fewer degrees of freedom than the corresponding tensor product elements, with a gap that increases as the degree of polynomial approximation order increases.
   \label{fig:DegsDofs}} 
  \Description{A DOF study for tensor product and trimmed serendipity elements.  We consider the case for all of $H^1$, \hdiv and \hcurl.}
\end{figure}

Accordingly, we have implemented the basis functions from \citet{gillette2019computational} using Firedrake's internal coding conventions and then carried out tests of various use cases in a reliable computational framework.
We show in \cref{fig:DegsDofs} how the number of degrees of freedom (DOFs) grow  in tensor product spaces ($\mathcal{Q}^-$) versus trimmed serendipity spaces ($\mathcal{S}^-$) in various element types, for increasing degree of polynomial approximation order.
Seeing how such quick back-of-the-envelope calculations might translate into significant computational savings requires the thorough implementation provided in this paper.


The main contributions of this paper are as follows:
\begin{enumerate}
\item We explain the implementation of trimmed serendipity elements within the Firedrake software environment and describe how users can easily employ the elements in their own code,
\item We validate the implementation of the trimmed serendipity finite elements by confirming they attain the theoretical bounds that are predicted in terms of convergence rates, and  
\item We examine the costs and benefits of using trimmed serendipity elements within the confines of typical test problems for numerical analysis by comparing them to tensor product elements.
\end{enumerate}

\section{Background and notation for trimmed serendipity elements}

    Among its many advantages, Finite Element Exterior Calculus (FEEC)~\cite{AFW2006,AFW2010,ABB2012} gives an easy, unified way to notate different element types.  The four best known families of elements are denoted $\calP^-_r \Lambda^k$, $\calP_r \Lambda^k$, $\calQ^-_r \Lambda^k$, and $\calS_r \Lambda^k$, which are, respectively, the trimmed polynomial, polynomial, tensor product, and serendipity elements of order $r$ using $k$-forms.  The $\calP$ and $\calP^-$ spaces are defined over meshes of simplices (triangles, tetrahedra, etc) while the $\calQ^-$ and $\calS$ spaces are defined over meshes of hypercubes (squares, cubes, etc). The optional notation $(\square_n)$ specifies that the space is constructed over the $n$-dimensional cube in $\R^n$, but we will frequently omit this addition if $n$ is clear from context.
    
    The mathematical results from FEEC regarding these four families synthesize decades of research into what constitutes a stable finite element, i.e.\ a numerical method that can be proven to converge to the correct solution of certain PDEs in certain norms at a prescribed rate, indicated by the subscript $r$.  In any dimension, the $0$-form spaces provide scalar-valued, $H^1$-conforming elements.  In 2D, $1$-forms can represent both \hcurl and \hdiv elements, depending on the orientation of the DOFs defined on the edges of a mesh.  In 3D, $1$-forms correspond to \hcurl elements while $2$-forms correspond to \hdiv elements.  Notably, the serendipity family $\calS_r\Lambda^k$ is the youngest, least implemented, and hence least understood among all these families.  In particular, the 1-form and 2-form (regular) serendipity elements in 3D were only characterized in 2014 by \citet{arnold2014finite}, whereas the equivalent tensor product elements were first described more than 30 years earlier in \citet{N1980,N1986}.

\begin{table}[htbp]
  \centering
  \caption{Trimmed serendipity elements on a reference cube in 3D, akin to the Periodic Table of the Finite Elements.  Columns show increasing form order from $k=0$ to $k=3$, corresponding to $H^1$,~\hcurl,~\hdiv, and $L^2$ conformity, respectively.  Rows show increasing order of approximation $r=1,2,3$. On the front face of each element, dots indicate the number of DOFs associated to each vertex, edge, or face of the cubical element. The number of DOFs associated to the interior of the cube is indicated with ``$\#$ int $=$.''  The total degree of freedom count for each element is shown to its right.\label{tab:tsfamily}}
\begin{tabular}{rm{.11\textwidth}m{.05\textwidth}m{.11\textwidth}m{.05\textwidth}m{.11\textwidth}m{.05\textwidth}m{.11\textwidth}m{.05\textwidth}}
 & $\qquad\calS_r^-\Lambda^0(\square_3)$ && $\qquad\calS_r^-\Lambda^1(\square_3)$  && $\qquad\calS_r^-\Lambda^2(\square_3)$  && $\qquad\calS_r^-\Lambda^3(\square_3)$ \\[2mm]
$r=1~~$ 
& \begin{tikzpicture}
\fill [cyan, fill opacity=0.1] (0,0,1) -- (1,0,1) -- (1,1,1) -- (0,1,1) -- cycle;
\fill [cyan, fill opacity=0.1] (0,0,0) -- (0,0,1) -- (1,0,1) -- (1,0,0) -- cycle;
\fill [cyan, fill opacity=0.1] (0,0,0) -- (1,0,0) -- (1,1,0) -- (0,1,0) -- cycle;
\fill [cyan, fill opacity=0.1] (0,0,0) -- (0,0,1) -- (0,1,1) -- (0,1,0) -- cycle;
\fill [cyan, fill opacity=0.1] (1,0,0) -- (1,0,1) -- (1,1,1) -- (1,1,0) -- cycle;
\fill [cyan, fill opacity=0.1] (0,1,0) -- (1,1,0) -- (1,1,1) -- (0,1,1) -- cycle;

\draw [line width=1] (0,1,0) -- (1,1,0); 
\draw [line width=1] (0,1,0) -- (0,1,1);
\draw [line width=1] (0,1,1) -- (1,1,1);
\draw [line width=1] (1,1,0) -- (1,1,1);
\draw [line width=1] (0,0,1) -- (0,1,1);
\draw [line width=1] (0,0,1) -- (1,0,1);
\draw [line width=1] (1,0,0) -- (1,0,1);
\draw [line width=1] (1,0,1) -- (1,1,1);
\draw [line width=1] (1,0,0) -- (1,1,0);
\draw [line width=1, dashed, opacity=0.25] (0,0,0) -- (1,0,0);
\draw [line width=1, dashed, opacity=0.25] (0,0,0) -- (0,1,0);
\draw [line width=1, dashed, opacity=0.25] (0,0,0) -- (0,0,1);

\filldraw [red] (1,0,0) circle (1.25pt);
\filldraw [red] (0,1,0) circle (1.25pt);
\filldraw [red] (0,0,1) circle (1.25pt);
\filldraw [red] (1,0,1) circle (1.25pt);
\filldraw [red] (1,1,1) circle (1.25pt);
\filldraw [red] (0,1,1) circle (1.25pt);
\filldraw [red] (1,1,0) circle (1.25pt);

\draw (0.6,0,1.5) node {\tiny{$\#$int = 0}};

\end{tikzpicture}
& \raisebox{8\height}{\Large 8}
& \begin{tikzpicture}
\fill [cyan, fill opacity=0.1] (0,0,1) -- (1,0,1) -- (1,1,1) -- (0,1,1) -- cycle;
\fill [cyan, fill opacity=0.1] (0,0,0) -- (0,0,1) -- (1,0,1) -- (1,0,0) -- cycle;
\fill [cyan, fill opacity=0.1] (0,0,0) -- (1,0,0) -- (1,1,0) -- (0,1,0) -- cycle;
\fill [cyan, fill opacity=0.1] (0,0,0) -- (0,0,1) -- (0,1,1) -- (0,1,0) -- cycle;
\fill [cyan, fill opacity=0.1] (1,0,0) -- (1,0,1) -- (1,1,1) -- (1,1,0) -- cycle;
\fill [cyan, fill opacity=0.1] (0,1,0) -- (1,1,0) -- (1,1,1) -- (0,1,1) -- cycle;

\draw [line width=1] (0,1,0) -- (1,1,0); 
\draw [line width=1] (0,1,0) -- (0,1,1);
\draw [line width=1] (0,1,1) -- (1,1,1);
\draw [line width=1] (1,1,0) -- (1,1,1);
\draw [line width=1] (0,0,1) -- (0,1,1);
\draw [line width=1] (0,0,1) -- (1,0,1);
\draw [line width=1] (1,0,0) -- (1,0,1);
\draw [line width=1] (1,0,1) -- (1,1,1);
\draw [line width=1] (1,0,0) -- (1,1,0);
\draw [line width=1, dashed, opacity=0.25] (0,0,0) -- (1,0,0);
\draw [line width=1, dashed, opacity=0.25] (0,0,0) -- (0,1,0);
\draw [line width=1, dashed, opacity=0.25] (0,0,0) -- (0,0,1);

\node[mark size=1.5pt,color=red, rotate=270] at (0.5,1,0) {\pgfuseplotmark{diamond*}};
\node[mark size=1.5pt,color=red, rotate=270] at (0.5,1,1) {\pgfuseplotmark{diamond*}};
\node[mark size=1.5pt,color=red, rotate=315] at (1,1,0.5) {\pgfuseplotmark{diamond*}};
\node[mark size=1.5pt,color=red, rotate=315] at (0,1,0.5) {\pgfuseplotmark{diamond*}};
\node[mark size=1.5pt,color=red, rotate=180] at (0,0.5,1) {\pgfuseplotmark{diamond*}};
\node[mark size=1.5pt,color=red, rotate=180] at (1,0.5,1) {\pgfuseplotmark{diamond*}};
\node[mark size=1.5pt,color=red, rotate=180] at (1,0.5,0) {\pgfuseplotmark{diamond*}};
\node[mark size=1.5pt,color=red, rotate=270] at (0.5,0,1) {\pgfuseplotmark{diamond*}};
\node[mark size=1.5pt,color=red, rotate=315] at (1,0,0.5) {\pgfuseplotmark{diamond*}};
\draw (0.6,0,1.5) node {\tiny{$\#$int = 0}};

\end{tikzpicture}
& \raisebox{8\height}{\Large 12}
& \begin{tikzpicture}
\fill [cyan, fill opacity=0.1] (0,0,1) -- (1,0,1) -- (1,1,1) -- (0,1,1) -- cycle;
\fill [cyan, fill opacity=0.1] (0,0,0) -- (0,0,1) -- (1,0,1) -- (1,0,0) -- cycle;
\fill [cyan, fill opacity=0.1] (0,0,0) -- (1,0,0) -- (1,1,0) -- (0,1,0) -- cycle;
\fill [cyan, fill opacity=0.1] (0,0,0) -- (0,0,1) -- (0,1,1) -- (0,1,0) -- cycle;
\fill [cyan, fill opacity=0.1] (1,0,0) -- (1,0,1) -- (1,1,1) -- (1,1,0) -- cycle;
\fill [cyan, fill opacity=0.1] (0,1,0) -- (1,1,0) -- (1,1,1) -- (0,1,1) -- cycle;

\draw [line width=1] (0,1,0) -- (1,1,0); 
\draw [line width=1] (0,1,0) -- (0,1,1);
\draw [line width=1] (0,1,1) -- (1,1,1);
\draw [line width=1] (1,1,0) -- (1,1,1);
\draw [line width=1] (0,0,1) -- (0,1,1);
\draw [line width=1] (0,0,1) -- (1,0,1);
\draw [line width=1] (1,0,0) -- (1,0,1);
\draw [line width=1] (1,0,1) -- (1,1,1);
\draw [line width=1] (1,0,0) -- (1,1,0);
\draw [line width=1, dashed, opacity=0.25] (0,0,0) -- (1,0,0);
\draw [line width=1, dashed, opacity=0.25] (0,0,0) -- (0,1,0);
\draw [line width=1, dashed, opacity=0.25] (0,0,0) -- (0,0,1);

\filldraw [red] (0.5,0.5,1) circle (1.25pt);
\filldraw [red] (1,0.5,0.5) circle (1.25pt);
\filldraw [red] (0.5,1,0.5) circle (1.25pt);

\draw (0.6,0,1.5) node {\tiny{$\#$int = 0}};

\end{tikzpicture}
& \raisebox{8\height}{\Large 6}
& \begin{tikzpicture}
\fill [cyan, fill opacity=0.1] (0,0,1) -- (1,0,1) -- (1,1,1) -- (0,1,1) -- cycle;
\fill [cyan, fill opacity=0.1] (0,0,0) -- (0,0,1) -- (1,0,1) -- (1,0,0) -- cycle;
\fill [cyan, fill opacity=0.1] (0,0,0) -- (1,0,0) -- (1,1,0) -- (0,1,0) -- cycle;
\fill [cyan, fill opacity=0.1] (0,0,0) -- (0,0,1) -- (0,1,1) -- (0,1,0) -- cycle;
\fill [cyan, fill opacity=0.1] (1,0,0) -- (1,0,1) -- (1,1,1) -- (1,1,0) -- cycle;
\fill [cyan, fill opacity=0.1] (0,1,0) -- (1,1,0) -- (1,1,1) -- (0,1,1) -- cycle;

\draw [line width=1] (0,1,0) -- (1,1,0); 
\draw [line width=1] (0,1,0) -- (0,1,1);
\draw [line width=1] (0,1,1) -- (1,1,1);
\draw [line width=1] (1,1,0) -- (1,1,1);
\draw [line width=1] (0,0,1) -- (0,1,1);
\draw [line width=1] (0,0,1) -- (1,0,1);
\draw [line width=1] (1,0,0) -- (1,0,1);
\draw [line width=1] (1,0,1) -- (1,1,1);
\draw [line width=1] (1,0,0) -- (1,1,0);
\draw [line width=1, dashed, opacity=0.25] (0,0,0) -- (1,0,0);
\draw [line width=1, dashed, opacity=0.25] (0,0,0) -- (0,1,0);
\draw [line width=1, dashed, opacity=0.25] (0,0,0) -- (0,0,1);

\draw (0.6,0,1.5) node {\tiny{$\#$int = 1}};

\end{tikzpicture}
& \raisebox{8\height}{\Large 1} \\
$r=2$
& \begin{tikzpicture}
\fill [cyan, fill opacity=0.1] (0,0,1) -- (1,0,1) -- (1,1,1) -- (0,1,1) -- cycle;
\fill [cyan, fill opacity=0.1] (0,0,0) -- (0,0,1) -- (1,0,1) -- (1,0,0) -- cycle;
\fill [cyan, fill opacity=0.1] (0,0,0) -- (1,0,0) -- (1,1,0) -- (0,1,0) -- cycle;
\fill [cyan, fill opacity=0.1] (0,0,0) -- (0,0,1) -- (0,1,1) -- (0,1,0) -- cycle;
\fill [cyan, fill opacity=0.1] (1,0,0) -- (1,0,1) -- (1,1,1) -- (1,1,0) -- cycle;
\fill [cyan, fill opacity=0.1] (0,1,0) -- (1,1,0) -- (1,1,1) -- (0,1,1) -- cycle;

\draw [line width=1] (0,1,0) -- (1,1,0); 
\draw [line width=1] (0,1,0) -- (0,1,1);
\draw [line width=1] (0,1,1) -- (1,1,1);
\draw [line width=1] (1,1,0) -- (1,1,1);
\draw [line width=1] (0,0,1) -- (0,1,1);
\draw [line width=1] (0,0,1) -- (1,0,1);
\draw [line width=1] (1,0,0) -- (1,0,1);
\draw [line width=1] (1,0,1) -- (1,1,1);
\draw [line width=1] (1,0,0) -- (1,1,0);
\draw [line width=1, dashed, opacity=0.25] (0,0,0) -- (1,0,0);
\draw [line width=1, dashed, opacity=0.25] (0,0,0) -- (0,1,0);
\draw [line width=1, dashed, opacity=0.25] (0,0,0) -- (0,0,1);

\filldraw [red] (1,0,0) circle (1.25pt);
\filldraw [red] (0,1,0) circle (1.25pt);
\filldraw [red] (0,0,1) circle (1.25pt);
\filldraw [red] (1,0,1) circle (1.25pt);
\filldraw [red] (1,1,1) circle (1.25pt);
\filldraw [red] (0,1,1) circle (1.25pt);
\filldraw [red] (1,1,0) circle (1.25pt);

\filldraw [red] (0,1,0.5) circle (1.25pt);
\filldraw [red] (0.5,1,1) circle (1.25pt);
\filldraw [red] (1,1,0.5) circle (1.25pt);
\filldraw [red] (0.5,1,0) circle (1.25pt);
\filldraw [red] (0,0.5,1) circle (1.25pt);
\filldraw [red] (1,0.5,1) circle (1.25pt);
\filldraw [red] (1,0.5,0) circle (1.25pt);
\filldraw [red] (0.5,0,1) circle (1.25pt);
\filldraw [red] (1,0,0.5) circle (1.25pt);

\draw (0.6,0,1.5) node {\tiny{$\#$int = 0}};

\end{tikzpicture}
& \raisebox{8\height}{\Large 20}
& \begin{tikzpicture}
\fill [cyan, fill opacity=0.1] (0,0,1) -- (1,0,1) -- (1,1,1) -- (0,1,1) -- cycle;
\fill [cyan, fill opacity=0.1] (0,0,0) -- (0,0,1) -- (1,0,1) -- (1,0,0) -- cycle;
\fill [cyan, fill opacity=0.1] (0,0,0) -- (1,0,0) -- (1,1,0) -- (0,1,0) -- cycle;
\fill [cyan, fill opacity=0.1] (0,0,0) -- (0,0,1) -- (0,1,1) -- (0,1,0) -- cycle;
\fill [cyan, fill opacity=0.1] (1,0,0) -- (1,0,1) -- (1,1,1) -- (1,1,0) -- cycle;
\fill [cyan, fill opacity=0.1] (0,1,0) -- (1,1,0) -- (1,1,1) -- (0,1,1) -- cycle;

\draw [line width=1] (0,1,0) -- (1,1,0); 
\draw [line width=1] (0,1,0) -- (0,1,1);
\draw [line width=1] (0,1,1) -- (1,1,1);
\draw [line width=1] (1,1,0) -- (1,1,1);
\draw [line width=1] (0,0,1) -- (0,1,1);
\draw [line width=1] (0,0,1) -- (1,0,1);
\draw [line width=1] (1,0,0) -- (1,0,1);
\draw [line width=1] (1,0,1) -- (1,1,1);
\draw [line width=1] (1,0,0) -- (1,1,0);
\draw [line width=1, dashed, opacity=0.25] (0,0,0) -- (1,0,0);
\draw [line width=1, dashed, opacity=0.25] (0,0,0) -- (0,1,0);
\draw [line width=1, dashed, opacity=0.25] (0,0,0) -- (0,0,1);

\node[mark size=1.5pt,color=red, rotate=270] at (0.33,1,1) {\pgfuseplotmark{diamond*}};
\node[mark size=1.5pt,color=red, rotate=270] at (0.66,1,1) {\pgfuseplotmark{diamond*}};
\node[mark size=1.5pt,color=red, rotate=270] at (0.33,0,1) {\pgfuseplotmark{diamond*}};
\node[mark size=1.5pt,color=red, rotate=270] at (0.66,0,1) {\pgfuseplotmark{diamond*}};
\node[mark size=1.5pt,color=red, rotate=180] at (0,0.33,1) {\pgfuseplotmark{diamond*}};
\node[mark size=1.5pt,color=red, rotate=180] at (0,0.66,1) {\pgfuseplotmark{diamond*}};
\node[mark size=1.5pt,color=red, rotate=180] at (1,0.33,1) {\pgfuseplotmark{diamond*}};
\node[mark size=1.5pt,color=red, rotate=180] at (1,0.66,1) {\pgfuseplotmark{diamond*}};
\node[mark size=1.5pt,color=red, rotate=135] at (1,0,0.66) {\pgfuseplotmark{diamond*}};
\node[mark size=1.5pt,color=red, rotate=135] at (1,0,0.33) {\pgfuseplotmark{diamond*}};
\node[mark size=1.5pt,color=red, rotate=135] at (1,1,0.66) {\pgfuseplotmark{diamond*}};
\node[mark size=1.5pt,color=red, rotate=135] at (1,1,0.33) {\pgfuseplotmark{diamond*}};
\node[mark size=1.5pt,color=red, rotate=135] at (0,1,0.66) {\pgfuseplotmark{diamond*}};
\node[mark size=1.5pt,color=red, rotate=135] at (0,1,0.33) {\pgfuseplotmark{diamond*}};
\node[mark size=1.5pt,color=red, rotate=90] at (0.33,1,0) {\pgfuseplotmark{diamond*}};
\node[mark size=1.5pt,color=red, rotate=90] at (0.66,1,0) {\pgfuseplotmark{diamond*}};
\node[mark size=1.5pt,color=red, rotate=180] at (1,0.33,0) {\pgfuseplotmark{diamond*}};
\node[mark size=1.5pt,color=red, rotate=180] at (1,0.66,0) {\pgfuseplotmark{diamond*}};

\filldraw [red] (0.33,1,0.5) circle (1.25pt);
\filldraw [red] (0.66,1,0.5) circle (1.25pt);
\filldraw [red] (0.66,0.5,1) circle (1.25pt);
\filldraw [red] (0.33,0.5,1) circle (1.25pt);
\filldraw [red] (1,0.5,0.33) circle (1.25pt);
\filldraw [red] (1,0.5,0.66) circle (1.25pt);

\draw (0.6,0,1.5) node {\tiny{$\#$int = 0}};

\end{tikzpicture}
& \raisebox{8\height}{\Large 36}
& \begin{tikzpicture}
\fill [cyan, fill opacity=0.1] (0,0,1) -- (1,0,1) -- (1,1,1) -- (0,1,1) -- cycle;
\fill [cyan, fill opacity=0.1] (0,0,0) -- (0,0,1) -- (1,0,1) -- (1,0,0) -- cycle;
\fill [cyan, fill opacity=0.1] (0,0,0) -- (1,0,0) -- (1,1,0) -- (0,1,0) -- cycle;
\fill [cyan, fill opacity=0.1] (0,0,0) -- (0,0,1) -- (0,1,1) -- (0,1,0) -- cycle;
\fill [cyan, fill opacity=0.1] (1,0,0) -- (1,0,1) -- (1,1,1) -- (1,1,0) -- cycle;
\fill [cyan, fill opacity=0.1] (0,1,0) -- (1,1,0) -- (1,1,1) -- (0,1,1) -- cycle;

\draw [line width=1] (0,1,0) -- (1,1,0); 
\draw [line width=1] (0,1,0) -- (0,1,1);
\draw [line width=1] (0,1,1) -- (1,1,1);
\draw [line width=1] (1,1,0) -- (1,1,1);
\draw [line width=1] (0,0,1) -- (0,1,1);
\draw [line width=1] (0,0,1) -- (1,0,1);
\draw [line width=1] (1,0,0) -- (1,0,1);
\draw [line width=1] (1,0,1) -- (1,1,1);
\draw [line width=1] (1,0,0) -- (1,1,0);
\draw [line width=1, dashed, opacity=0.25] (0,0,0) -- (1,0,0);
\draw [line width=1, dashed, opacity=0.25] (0,0,0) -- (0,1,0);
\draw [line width=1, dashed, opacity=0.25] (0,0,0) -- (0,0,1);

\filldraw [red] (0.4,0.4,1) circle (1.25pt);
\filldraw [red] (0.6,0.4,1) circle (1.25pt);
\filldraw [red] (0.5,0.6,1) circle (1.25pt);

\filldraw [red] (1,0.4,0.4) circle (1.25pt);
\filldraw [red] (1,0.4,0.6) circle (1.25pt);
\filldraw [red] (1,0.6,0.5) circle (1.25pt);

\filldraw [red] (0.3,1,0.4) circle (1.25pt);
\filldraw [red] (0.6,1,0.4) circle (1.25pt);
\filldraw [red] (0.5,1,0.6) circle (1.25pt);

\draw (0.6,0,1.5) node {\tiny{$\#$int = 3}};

\end{tikzpicture}
& \raisebox{8\height}{\Large 21}
& \begin{tikzpicture}
\fill [cyan, fill opacity=0.1] (0,0,1) -- (1,0,1) -- (1,1,1) -- (0,1,1) -- cycle;
\fill [cyan, fill opacity=0.1] (0,0,0) -- (0,0,1) -- (1,0,1) -- (1,0,0) -- cycle;
\fill [cyan, fill opacity=0.1] (0,0,0) -- (1,0,0) -- (1,1,0) -- (0,1,0) -- cycle;
\fill [cyan, fill opacity=0.1] (0,0,0) -- (0,0,1) -- (0,1,1) -- (0,1,0) -- cycle;
\fill [cyan, fill opacity=0.1] (1,0,0) -- (1,0,1) -- (1,1,1) -- (1,1,0) -- cycle;
\fill [cyan, fill opacity=0.1] (0,1,0) -- (1,1,0) -- (1,1,1) -- (0,1,1) -- cycle;

\draw [line width=1] (0,1,0) -- (1,1,0); 
\draw [line width=1] (0,1,0) -- (0,1,1);
\draw [line width=1] (0,1,1) -- (1,1,1);
\draw [line width=1] (1,1,0) -- (1,1,1);
\draw [line width=1] (0,0,1) -- (0,1,1);
\draw [line width=1] (0,0,1) -- (1,0,1);
\draw [line width=1] (1,0,0) -- (1,0,1);
\draw [line width=1] (1,0,1) -- (1,1,1);
\draw [line width=1] (1,0,0) -- (1,1,0);
\draw [line width=1, dashed, opacity=0.25] (0,0,0) -- (1,0,0);
\draw [line width=1, dashed, opacity=0.25] (0,0,0) -- (0,1,0);
\draw [line width=1, dashed, opacity=0.25] (0,0,0) -- (0,0,1);

\draw (0.6,0,1.5) node {\tiny{$\#$int = 4}};

\end{tikzpicture}
& \raisebox{8\height}{\Large 4} \\
$r=3$
& \begin{tikzpicture}
\fill [cyan, fill opacity=0.1] (0,0,1) -- (1,0,1) -- (1,1,1) -- (0,1,1) -- cycle;
\fill [cyan, fill opacity=0.1] (0,0,0) -- (0,0,1) -- (1,0,1) -- (1,0,0) -- cycle;
\fill [cyan, fill opacity=0.1] (0,0,0) -- (1,0,0) -- (1,1,0) -- (0,1,0) -- cycle;
\fill [cyan, fill opacity=0.1] (0,0,0) -- (0,0,1) -- (0,1,1) -- (0,1,0) -- cycle;
\fill [cyan, fill opacity=0.1] (1,0,0) -- (1,0,1) -- (1,1,1) -- (1,1,0) -- cycle;
\fill [cyan, fill opacity=0.1] (0,1,0) -- (1,1,0) -- (1,1,1) -- (0,1,1) -- cycle;

\draw [line width=1] (0,1,0) -- (1,1,0); 
\draw [line width=1] (0,1,0) -- (0,1,1);
\draw [line width=1] (0,1,1) -- (1,1,1);
\draw [line width=1] (1,1,0) -- (1,1,1);
\draw [line width=1] (0,0,1) -- (0,1,1);
\draw [line width=1] (0,0,1) -- (1,0,1);
\draw [line width=1] (1,0,0) -- (1,0,1);
\draw [line width=1] (1,0,1) -- (1,1,1);
\draw [line width=1] (1,0,0) -- (1,1,0);
\draw [line width=1, dashed, opacity=0.25] (0,0,0) -- (1,0,0);
\draw [line width=1, dashed, opacity=0.25] (0,0,0) -- (0,1,0);
\draw [line width=1, dashed, opacity=0.25] (0,0,0) -- (0,0,1);

\filldraw [red] (1,0,0) circle (1.25pt);
\filldraw [red] (0,1,0) circle (1.25pt);
\filldraw [red] (0,0,1) circle (1.25pt);
\filldraw [red] (1,0,1) circle (1.25pt);
\filldraw [red] (1,1,1) circle (1.25pt);
\filldraw [red] (0,1,1) circle (1.25pt);
\filldraw [red] (1,1,0) circle (1.25pt);

\filldraw [red] (0,1,0.33) circle (1.25pt);
\filldraw [red] (0.33,1,1) circle (1.25pt);
\filldraw [red] (1,1,0.33) circle (1.25pt);
\filldraw [red] (0.33,1,0) circle (1.25pt);
\filldraw [red] (0,0.33,1) circle (1.25pt);
\filldraw [red] (1,0.33,1) circle (1.25pt);
\filldraw [red] (1,0.33,0) circle (1.25pt);
\filldraw [red] (0.33,0,1) circle (1.25pt);
\filldraw [red] (1,0,0.33) circle (1.25pt);

\filldraw [red] (0,1,0.66) circle (1.25pt);
\filldraw [red] (0.66,1,1) circle (1.25pt);
\filldraw [red] (1,1,0.66) circle (1.25pt);
\filldraw [red] (0.66,1,0) circle (1.25pt);
\filldraw [red] (0,0.66,1) circle (1.25pt);
\filldraw [red] (1,0.66,1) circle (1.25pt);
\filldraw [red] (1,0.66,0) circle (1.25pt);
\filldraw [red] (0.66,0,1) circle (1.25pt);
\filldraw [red] (1,0,0.66) circle (1.25pt);

\draw (0.6,0,1.5) node {\tiny{$\#$int = 0}};

\end{tikzpicture}
& \raisebox{8\height}{\Large 32}
& \begin{tikzpicture}
\fill [cyan, fill opacity=0.1] (0,0,1) -- (1,0,1) -- (1,1,1) -- (0,1,1) -- cycle;
\fill [cyan, fill opacity=0.1] (0,0,0) -- (0,0,1) -- (1,0,1) -- (1,0,0) -- cycle;
\fill [cyan, fill opacity=0.1] (0,0,0) -- (1,0,0) -- (1,1,0) -- (0,1,0) -- cycle;
\fill [cyan, fill opacity=0.1] (0,0,0) -- (0,0,1) -- (0,1,1) -- (0,1,0) -- cycle;
\fill [cyan, fill opacity=0.1] (1,0,0) -- (1,0,1) -- (1,1,1) -- (1,1,0) -- cycle;
\fill [cyan, fill opacity=0.1] (0,1,0) -- (1,1,0) -- (1,1,1) -- (0,1,1) -- cycle;

\draw [line width=1] (0,1,0) -- (1,1,0); 
\draw [line width=1] (0,1,0) -- (0,1,1);
\draw [line width=1] (0,1,1) -- (1,1,1);
\draw [line width=1] (1,1,0) -- (1,1,1);
\draw [line width=1] (0,0,1) -- (0,1,1);
\draw [line width=1] (0,0,1) -- (1,0,1);
\draw [line width=1] (1,0,0) -- (1,0,1);
\draw [line width=1] (1,0,1) -- (1,1,1);
\draw [line width=1] (1,0,0) -- (1,1,0);
\draw [line width=1, dashed, opacity=0.25] (0,0,0) -- (1,0,0);
\draw [line width=1, dashed, opacity=0.25] (0,0,0) -- (0,1,0);
\draw [line width=1, dashed, opacity=0.25] (0,0,0) -- (0,0,1);

\node[mark size=1.5pt,color=red, rotate=270] at (0.3,1,1) {\pgfuseplotmark{diamond*}};
\node[mark size=1.5pt,color=red, rotate=270] at (0.5,1,1) {\pgfuseplotmark{diamond*}};
\node[mark size=1.5pt,color=red, rotate=270] at (0.7,1,1) {\pgfuseplotmark{diamond*}};

\node[mark size=1.5pt,color=red, rotate=270] at (0.3,0,1) {\pgfuseplotmark{diamond*}};
\node[mark size=1.5pt,color=red, rotate=270] at (0.5,0,1) {\pgfuseplotmark{diamond*}};
\node[mark size=1.5pt,color=red, rotate=270] at (0.7,0,1) {\pgfuseplotmark{diamond*}};

\node[mark size=1.5pt,color=red, rotate=180] at (0,0.3,1) {\pgfuseplotmark{diamond*}};
\node[mark size=1.5pt,color=red, rotate=180] at (0,0.5,1) {\pgfuseplotmark{diamond*}};
\node[mark size=1.5pt,color=red, rotate=180] at (0,0.7,1) {\pgfuseplotmark{diamond*}};

\node[mark size=1.5pt,color=red, rotate=180] at (1,0.3,1) {\pgfuseplotmark{diamond*}};
\node[mark size=1.5pt,color=red, rotate=180] at (1,0.5,1) {\pgfuseplotmark{diamond*}};
\node[mark size=1.5pt,color=red, rotate=180] at (1,0.7,1) {\pgfuseplotmark{diamond*}};

\node[mark size=1.5pt,color=red, rotate=135] at (1,0,0.3) {\pgfuseplotmark{diamond*}};
\node[mark size=1.5pt,color=red, rotate=135] at (1,0,0.5) {\pgfuseplotmark{diamond*}};
\node[mark size=1.5pt,color=red, rotate=135] at (1,0,0.7) {\pgfuseplotmark{diamond*}};

\node[mark size=1.5pt,color=red, rotate=135] at (1,1,0.3) {\pgfuseplotmark{diamond*}};
\node[mark size=1.5pt,color=red, rotate=135] at (1,1,0.5) {\pgfuseplotmark{diamond*}};
\node[mark size=1.5pt,color=red, rotate=135] at (1,1,0.7) {\pgfuseplotmark{diamond*}};

\node[mark size=1.5pt,color=red, rotate=135] at (0,1,0.3) {\pgfuseplotmark{diamond*}};
\node[mark size=1.5pt,color=red, rotate=135] at (0,1,0.5) {\pgfuseplotmark{diamond*}};
\node[mark size=1.5pt,color=red, rotate=135] at (0,1,0.7) {\pgfuseplotmark{diamond*}};

\node[mark size=1.5pt,color=red, rotate=90] at (0.3,1,0) {\pgfuseplotmark{diamond*}};
\node[mark size=1.5pt,color=red, rotate=90] at (0.5,1,0) {\pgfuseplotmark{diamond*}};
\node[mark size=1.5pt,color=red, rotate=90] at (0.7,1,0) {\pgfuseplotmark{diamond*}};

\node[mark size=1.5pt,color=red, rotate=180] at (1,0.3,0) {\pgfuseplotmark{diamond*}};
\node[mark size=1.5pt,color=red, rotate=180] at (1,0.5,0) {\pgfuseplotmark{diamond*}};
\node[mark size=1.5pt,color=red, rotate=180] at (1,0.7,0) {\pgfuseplotmark{diamond*}};

\filldraw [red] (0.5,1,0.5) circle (1.25pt);
\filldraw [red] (0.66,1,0.66) circle (1.25pt);
\filldraw [red] (0.66,1,0.33) circle (1.25pt);
\filldraw [red] (0.33,1,0.66) circle (1.25pt);
\filldraw [red] (0.33,1,0.33) circle (1.25pt);

\filldraw [red] (0.5,0.5,1) circle (1.25pt);
\filldraw [red] (0.33,0.33,1) circle (1.25pt);
\filldraw [red] (0.33,0.66,1) circle (1.25pt);
\filldraw [red] (0.66,0.33,1) circle (1.25pt);
\filldraw [red] (0.66,0.66,1) circle (1.25pt);

\filldraw [red] (1,0.5,0.5) circle (1.25pt);
\filldraw [red] (1,0.33,0.66) circle (1.25pt);
\filldraw [red] (1,0.33,0.33) circle (1.25pt);
\filldraw [red] (1,0.66,0.66) circle (1.25pt);
\filldraw [red] (1,0.66,0.33) circle (1.25pt);

\draw (0.6,0,1.5) node {\tiny{$\#$int = 0}};

\end{tikzpicture}
& \raisebox{8\height}{\Large 66}
& \begin{tikzpicture}
\fill [cyan, fill opacity=0.1] (0,0,1) -- (1,0,1) -- (1,1,1) -- (0,1,1) -- cycle;
\fill [cyan, fill opacity=0.1] (0,0,0) -- (0,0,1) -- (1,0,1) -- (1,0,0) -- cycle;
\fill [cyan, fill opacity=0.1] (0,0,0) -- (1,0,0) -- (1,1,0) -- (0,1,0) -- cycle;
\fill [cyan, fill opacity=0.1] (0,0,0) -- (0,0,1) -- (0,1,1) -- (0,1,0) -- cycle;
\fill [cyan, fill opacity=0.1] (1,0,0) -- (1,0,1) -- (1,1,1) -- (1,1,0) -- cycle;
\fill [cyan, fill opacity=0.1] (0,1,0) -- (1,1,0) -- (1,1,1) -- (0,1,1) -- cycle;

\draw [line width=1] (0,1,0) -- (1,1,0); 
\draw [line width=1] (0,1,0) -- (0,1,1);
\draw [line width=1] (0,1,1) -- (1,1,1);
\draw [line width=1] (1,1,0) -- (1,1,1);
\draw [line width=1] (0,0,1) -- (0,1,1);
\draw [line width=1] (0,0,1) -- (1,0,1);
\draw [line width=1] (1,0,0) -- (1,0,1);
\draw [line width=1] (1,0,1) -- (1,1,1);
\draw [line width=1] (1,0,0) -- (1,1,0);
\draw [line width=1, dashed, opacity=0.25] (0,0,0) -- (1,0,0);
\draw [line width=1, dashed, opacity=0.25] (0,0,0) -- (0,1,0);
\draw [line width=1, dashed, opacity=0.25] (0,0,0) -- (0,0,1);

\filldraw [red] (0.4,0.5,1) circle (1.25pt);
\filldraw [red] (0.6,0.5,1) circle (1.25pt);
\filldraw [red] (0.5,0.7,1) circle (1.25pt);
\filldraw [red] (0.3,0.3,1) circle (1.25pt);
\filldraw [red] (0.5,0.3,1) circle (1.25pt);
\filldraw [red] (0.7,0.3,1) circle (1.25pt);

\filldraw [red] (1,0.5,0.4) circle (1.25pt);
\filldraw [red] (1,0.5,0.6) circle (1.25pt);
\filldraw [red] (1,0.7,0.5) circle (1.25pt);
\filldraw [red] (1,0.3,0.3) circle (1.25pt);
\filldraw [red] (1,0.3,0.5) circle (1.25pt);
\filldraw [red] (1,0.3,0.7) circle (1.25pt);

\filldraw [red] (0.4,1,0.5) circle (1.25pt);
\filldraw [red] (0.6,1,0.5) circle (1.25pt);
\filldraw [red] (0.5,1,0.7) circle (1.25pt);
\filldraw [red] (0.3,1,0.3) circle (1.25pt);
\filldraw [red] (0.5,1,0.3) circle (1.25pt);
\filldraw [red] (0.7,1,0.3) circle (1.25pt);

\draw (0.6,0,1.5) node {\tiny{$\#$int = 9}};

\end{tikzpicture}
& \raisebox{8\height}{\Large 45}
& \begin{tikzpicture}
\fill [cyan, fill opacity=0.1] (0,0,1) -- (1,0,1) -- (1,1,1) -- (0,1,1) -- cycle;
\fill [cyan, fill opacity=0.1] (0,0,0) -- (0,0,1) -- (1,0,1) -- (1,0,0) -- cycle;
\fill [cyan, fill opacity=0.1] (0,0,0) -- (1,0,0) -- (1,1,0) -- (0,1,0) -- cycle;
\fill [cyan, fill opacity=0.1] (0,0,0) -- (0,0,1) -- (0,1,1) -- (0,1,0) -- cycle;
\fill [cyan, fill opacity=0.1] (1,0,0) -- (1,0,1) -- (1,1,1) -- (1,1,0) -- cycle;
\fill [cyan, fill opacity=0.1] (0,1,0) -- (1,1,0) -- (1,1,1) -- (0,1,1) -- cycle;

\draw [line width=1] (0,1,0) -- (1,1,0); 
\draw [line width=1] (0,1,0) -- (0,1,1);
\draw [line width=1] (0,1,1) -- (1,1,1);
\draw [line width=1] (1,1,0) -- (1,1,1);
\draw [line width=1] (0,0,1) -- (0,1,1);
\draw [line width=1] (0,0,1) -- (1,0,1);
\draw [line width=1] (1,0,0) -- (1,0,1);
\draw [line width=1] (1,0,1) -- (1,1,1);
\draw [line width=1] (1,0,0) -- (1,1,0);
\draw [line width=1, dashed, opacity=0.25] (0,0,0) -- (1,0,0);
\draw [line width=1, dashed, opacity=0.25] (0,0,0) -- (0,1,0);
\draw [line width=1, dashed, opacity=0.25] (0,0,0) -- (0,0,1);

\draw (0.6,0,1.5) node {\tiny{$\#$int = 10}};

\end{tikzpicture}
& \raisebox{8\height}{\Large 10}
\end{tabular}
\end{table}

\subsection{Trimmed Serendipity Elements}
    The trimmed serendipity family is an even newer addition to this collection that attains a key optimality condition arising from the FEEC framework.
	\citet{christiansen2016constructions} computed the minimal possible dimensions for an exact sequence of conforming finite element spaces on cubes that contained $\calP_r\Lambda^k$ for each $k$.
	\citet{gillette2019trimmed} identified polynomial differential form spaces with these prescribed dimensions and approximation power, denoting them trimmed serendipity elements with the notation $S^-_r\Lambda^k$.  
	Thus, trimmed serendipity elements represent the cheapest possible way to get an order $r$ conforming finite element method, if cost is only measured in terms of number of degrees of freedom.

In order to test if this benefit translated to computational speedups, Gillette, Kloefkorn and Sanders gave a systematic way to build the computational basis for each of these elements for dimensions $n = 2, 3$, $k=0, 1, 2, 3$-forms, and arbitrary order $r \geq 1$~\cite{gillette2019computational}.
These ``computational bases'' are well-suited for implementation since each basis function is associated to a unique mesh identity - i.e.\ a specific vertex, edge, face (for cubes) or element interior.
The required geometric localization of DOFs is visualized for low orders in 3D in \cref{tab:tsfamily}, arranged equivalently to the Periodic Table of the Finite Elements.
We note that neither the particular bases defined in \citet{gillette2019computational} nor any other general implementation of trimmed serendipity elements has been attempted prior to this paper.

  \subsubsection{Scalar Trimmed Serendipity Elements}
  The scalar-valued trimmed serendipity elements that are represented by $0$-forms are used as the shape functions for an $H^1$-conforming finite element space.  These are denoted by $\calS_r^-\Lambda^0$ and are identical to the scalar-valued serendipity elements from the Periodic Table of Finite Elements, i.e.\ $\calS_r^-\Lambda^0(\square_n) = \calS_r\Lambda^0(\square_n)$ for any $n$.  Arnold and Awanou provided a simple description of the functions in $\mathcal{S}_r\Lambda^0$ as the span of all monomials of ``superlinear degree'' less than or equal to $r$~\cite{arnold2011serendipity}. 
  
  Likewise, the scalar-valued trimmed serendipity elements that are represented by $n$-forms create $L^2$-conforming finite element spaces.  These are denoted by $\calS_r^-\Lambda^n(\square_n)$, and here we have the equality $\calS_r^-\Lambda^n(\square_n) = \calS_{r+1}\Lambda^n(\square_n)$.  In terms of the Periodic Table of Finite Elements, these are the \textbf{dPc}$_r$ spaces.  The shape functions for these spaces are simply the space of order $r$ polynomials.  Since no inter-element continuity is needed for $L^2$-conformity, we have the additional equivalence $\calS_r^-\Lambda^n(\square_n) = \calP_{r+1}\Lambda^n(\square_n)$.
  
  \subsubsection{Vector-valued Trimmed Serendipity Elements}

	The trimmed serendipity elements are truly distinct from regular serendipity spaces for $k$ values $0<k<n$.  Here, we will only consider dimensions $n=2$ and $n=3$, where the $k$-form spaces can be identified as vector-valued finite elements.  	
	In 2D, the vector-valued spaces $\calS_r^-\Lambda^1(\square_2)$ bear close relation to the Arbogast-Correa elements~\citep{arbogast2016two}, as explained in \citet[Prop 2.2]{gillette2019trimmed}.
	In 3D, the vector valued spaces $\calS_r^-\Lambda^1(\square_3)$ and $\calS_r^-\Lambda^2(\square_3)$ were also characterized by \citet{CF2016}, as explained in \citet[Prop 2.3]{gillette2019trimmed}.

	A major reason that the vector trimmed serendipity elements have only recently been considered in the mathematical literature is that their DOF per element count is complicated.
	As evidenced by \cref{tab:tsfamily}, \textit{certain} DOFs grow in predictable patterns with $r$.
	For instance, elements in the 1-form family, $\calS_r^-\Lambda^1(\square_3)$, have exactly $r$ DOFs per edge of the cube (corresponding to ``order $r$'' approximation on edges) and elements in the 2-form family, $\calS_r^-\Lambda^2(\square_3)$, have ${\displaystyle {r+1}\choose 2}$ DOFs per face (corresponding to ``order $r$'' approximation on faces).
	However, the 2-form family has the more obscure quantity of $(r^3 - 2r^2 + 3r)/2$ DOFs associated to the interior of an element (for $r>1$).
	This growth pattern is recognized as sequence A064808 by the On-line Encyclopedia of Integer Sequences~\cite{Sloane2018} and is in agreement with the general formulae presented in \citet{gillette2019trimmed}, but a natural geometric interpretation remains elusive.
	Equally unexpected patterns are evident in the growth of DOFs on faces and interiors of the $\calS_r^-\Lambda^1(\square_3)$ family.
	As we will discuss in the next section, Firedrake makes the creation and incorporation of such unusual DOF growth patterns simple for the developer, and thus opens these elements to numerical testing for the first time.

  \section{Building capacity for serendipity element types in Firedrake}
  \label{sec:buildcap}

  Firedrake uses FIAT \cite{kirby2004algorithm,kirby2012fiat} to provide finite element basis functions on reference elements.  
  To implement a new element in FIAT, we must provide both rules to tabulate basis functions and their derivatives at reference element points and a data structure that assigns basis functions to particular reference element facets.
  Said element is then made available in Firedrake by providing a symbolic name (in UFL \cite{alnaes2014unified}) and a translation from symbolic name to concrete implementation in the form compiler TSFC \cite{homolya2018tsfc}.  While FIAT initially considered a very wide range of finite elements~\cite{kirby2012common}, it would seek to express their bases as linear combinations of orthogonal polynomials.  However, for some elements, it is easier to directly describe the basis functions.
  We follow the construction of \citet{gillette2019computational} which provides explicit formulae for the basis functions for each of the trimmed serendipity elements and directly implement tabulation of the basis functions. To provide tabulations of derivatives, we implement the basis functions symbolically with SymPy \cite{sympy2017} and compute derivatives symbolically.

We use the decompositions from \citet{gillette2019computational} to group the basis functions according to the geometric portions of a reference mesh element (vertices, edges, faces, or cell interiors).  For instance, a basis for $\calS^-_r\Lambda^1(\square_3)$ -- the trimmed serendipity \hcurl-conforming element in 3D -- can be decomposed as
   \begin{equation}\label{eq:HCurlTrimmedSerendipityBasis}
   \calS^-_r\Lambda^1(\square_3) =   \underbrace{\left[\bigoplus_{i=0}^{r-1} E_i \Lambda^1\right]}_{\text{edge functions}}\oplus\underbrace{\left[ \bigoplus_{i=2}^{r-1}F_i \Lambda^1\right] \oplus \left[\tilde{F}_r \Lambda^1\right]}_{\text{face functions}}\oplus\underbrace{\left[ \bigoplus_{i=4}^{r-1}I_i \Lambda^1 \right] \oplus \left[\tilde{I}_r \Lambda^1\right]}_{\text{interior functions}}.
   \end{equation}
  Subsets in these decompositions denoted with an $E$, $F$, or $I$ are defined on edges, faces, and interior, respectively, of the cubical cell in 3D.  The subsets $\tilde{F}$ and $\tilde{I}$ are extra sets of basis functions defined on the faces or the interior that follow a different pattern in their definitions than those of the functions in $E$, $F$, or $I$. 

To see how this plays out in the Firedrake implementation, consider the \hcurl elements for trimmed serendipity at order $r=2$ in $n=3$ dimensions, indicating the space $S_2^- \Lambda^1$ in 3D.
The space $\tilde{I}_r \Lambda^1$ is defined to be empty for $r<4$, so there are only two sets of functions to include in this case: one set associated to edges of the reference element (the $E_i\Lambda^1$ sum) and one set for the faces of the reference element (the $\tilde{F}_2\Lambda^1$ functions).
According to \cref{eq:HCurlTrimmedSerendipityBasis}, the basis functions can be decomposed as
\begin{equation}\label{eq:HcurlRtwo}
   \calS^-_2\Lambda^1(\square_3) =    \left[\bigoplus_{i=0}^{1} E_i \Lambda^1\right] \oplus \left[\tilde{F}_2 \Lambda^1\right]\oplus \left[\tilde{I}_2 \Lambda^1\right].
   \end{equation}

We then implement these in Firedrake as follows.  The first step is to determine the number of DOFs we will need on the reference element where we will define the basis functions.  To do this, we count the DOFs for each mesh entity on the reference element (vertex, edge, face, and interior).  For example, $\mathcal{S}_2^- \Lambda^1(\square_3)$ should have no DOFs at the vertices (these only come into play for the $H^1$-conforming elements), two DOFs on each edge, two DOFs on each face, and no DOFs on the interior.  This agrees with \cref{eq:HcurlRtwo}, where $E_i\Lambda^1$ supplies one DOF to each edge for each $i=0, 1$, and then $\tilde{F}_2\Lambda^1$ supplies two DOFs to each face.  An illustration of this can be seen in the second column, second row of \cref{tab:tsfamily}.

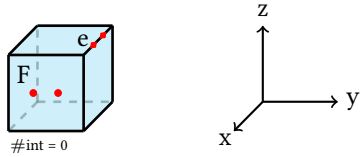
\begin{figure}[htbp]
\begin{tikzpicture}
\fill [cyan, fill opacity=0.1] (0,0,1) -- (1,0,1) -- (1,1,1) -- (0,1,1) -- cycle;
\fill [cyan, fill opacity=0.1] (0,0,0) -- (0,0,1) -- (1,0,1) -- (1,0,0) -- cycle;
\fill [cyan, fill opacity=0.1] (0,0,0) -- (1,0,0) -- (1,1,0) -- (0,1,0) -- cycle;
\fill [cyan, fill opacity=0.1] (0,0,0) -- (0,0,1) -- (0,1,1) -- (0,1,0) -- cycle;
\fill [cyan, fill opacity=0.1] (1,0,0) -- (1,0,1) -- (1,1,1) -- (1,1,0) -- cycle;
\fill [cyan, fill opacity=0.1] (0,1,0) -- (1,1,0) -- (1,1,1) -- (0,1,1) -- cycle;

\draw [line width=1] (0,1,0) -- (1,1,0); 
\draw [line width=1] (0,1,0) -- (0,1,1);
\draw [line width=1] (0,1,1) -- (1,1,1);
\draw [line width=1] (1,1,0) -- (1,1,1);
\draw [line width=1] (0,0,1) -- (0,1,1);
\draw [line width=1] (0,0,1) -- (1,0,1);
\draw [line width=1] (1,0,0) -- (1,0,1);
\draw [line width=1] (1,0,1) -- (1,1,1);
\draw [line width=1] (1,0,0) -- (1,1,0);
\draw [line width=1, dashed, opacity=0.25] (0,0,0) -- (1,0,0);
\draw [line width=1, dashed, opacity=0.25] (0,0,0) -- (0,1,0);
\draw [line width=1, dashed, opacity=0.25] (0,0,0) -- (0,0,1);

\draw (0.8,1,0.5) node {e};
\draw (0.2,0.75,1) node {F};

\filldraw [red] (0.66,0.5,1) circle (1.25pt);
\filldraw [red] (0.33,0.5,1) circle (1.25pt);

\draw[thick, ->] (3,0,0) -- (3,0,1);
\draw[thick, ->] (3,0,0) -- (3,1,0);
\draw[thick, ->] (3,0,0) -- (4,0,0);
\draw (3,0,1.3) node {x};
\draw (3,1.2,0) node {z};
\draw (4.2,0,0) node {y};


\node[mark size=1.5pt,color=red, rotate=135] at (1,1,0.66) {\pgfuseplotmark{diamond*}};
\node[mark size=1.5pt,color=red, rotate=135] at (1,1,0.33) {\pgfuseplotmark{diamond*}};




\draw (0.6,0,1.5) node {\tiny{$\#$int = 0}};

\end{tikzpicture}

\caption{The reference cube $[0,1]^3$ is shown, with coordinate axes indicated.  The edge $e$ lies in the intersection of the planes $y=1$ and $z=1$.  To ensure proper continuity, basis functions associated to $e$ must vanish on all edges of the cube \textit{except} $e$.  Examples of such functions for $e$ are described in detail in the text.  Additional examples of functions associated to the face $F$ are also given.\label{fig:ReferenceCube}}
\Description{An illustration of a reference element for trimmed serendipity 1-forms.  The dots indicate DOFs on the faces or edges.}
\end{figure}

With the number of DOFs assigned to each mesh entity in the reference element, we can then define the basis functions.  The order of definition is important so that basis functions are matched to the proper mesh entity.  Note that \cref{eq:HcurlRtwo} doesn't explicitly give a way to order the basis functions.  Instead, we need to use the properties of the basis functions to determine the correct ordering.   This is best illustrated by an example. Two of the ``edge'' basis functions that are contained in the sum of the $E_i\Lambda^1$ sets are $(y+1)(z+1)dx$ and $x(y+1)(z+1)dx$.  Notice that these functions have no $dy$ or $dz$ portions.  Therefore, these function vanish on any edge \textit{not} parallel to the $x$ axis.  Further, the polynomial coefficients of these forms indicate that they also vanish on the planes $y=-1$ and $z=-1$.  Thus, the only edge of the cube on which these functions do \textit{not} vanish is the edge contained in the line $\{y=1\}\cap\{z=1\}$.  This edge is shown in blue and labeled with an ``e'' in \cref{fig:ReferenceCube}.

This process is what determines the ordering for the basis functions.  If, in the first step described above where we are assigning DOFs to mesh entities, we assign the first two DOFs to be on the edge of the reference element where $y = 1$ and $z = 1$, then we must define the basis functions $(y+1)(z+1)dx$ and $x(y+1)(z+1)dx$ as our first two basis functions.

While the formulas above are taken from \citet{gillette2019computational}, these monomials have poor conditioning at higher orders.
Therefore, we use Legendre polynomials obtained symbolically from SymPy via the \texttt{legendre} function denoted \texttt{leg}.  Then we write the differential forms in vector notation.  For our examples of $(y+1)(z+1)dx$ and $x(y+1)(z+1)dx$, we get
\texttt{tuple([(leg(j,x\_mid)*dz[1]*dy[1],0,0)])}, where \texttt{leg(j, x\_mid)} is computed at the midpoint between two vertices and used for the $1$ or $x$ coefficient, and the \texttt{dy[1]} and \texttt{dz[1]} are used for the values of $(y+1)$ and $(z+1)$, respectively.  After repeating this process for each of the edges and associated basis functions, we then also do a similar process for the face functions in the set $\tilde{F}_2 \Lambda^1$.  

This process changes only slightly if we instead had considered the the $2$-forms $\mathcal{S}^-_2 \Lambda^2(\square_3)$.  In this scenario, we would have the basis functions given by the equation

\begin{equation*}
\calS^-_2\Lambda^2(\square_3) =    \left[\bigoplus_{i=0}^{1} F_i \Lambda^2\right] \oplus \left[\tilde{I}_2 \Lambda^2\right].
\end{equation*}

\noindent In this case, one of the basis functions is $y^jz^k(x+1)dydz$.  The $dydz$ represents a $2$-form, which vanishes on any face \textit{not} parallel to the $yz$ plane.  This leaves only the faces contained in the planes $x=1$ or $x=-1$ as possibilities for association. 
As in the 1-form example above, the polynomial coefficient of the form indicates that this function will vanish on an additional mesh entity, namely, the face contained in the plane $x=-1$, in this case.
Hence, we associate this function with the face contained in $x=1$ (labeled with an ``F'' in \cref{fig:ReferenceCube}).  The rest of the process for defining the $2$-form basis functions is similar to the process for the $1$-form basis functions.

\begin{table}[htbp]
  \centering
  \caption{A translation between FEEC and Firedrake usage names for tensor product and trimmed serendipity elements.  In 2D, the $0$-forms are $H^1$ conforming spaces, the $1$-forms are \hcurl and \hdiv conforming spaces (dependent upon oreintation of the DOFs), and the $2$-forms are $L^2$ conforming spaces.  For 3D, the $0$-forms are $H^1$ conforming spaces, the $1$-forms are \hcurl conforming spaces, $2$-forms are \hdiv conforming spaces, and $3$-forms are $L^2$ conforming spaces.\label{tab:FiredrakeNames}}
\begin{tabular}{llll}
  \toprule
   FEEC & UFL name (2D) & UFL name (3D) \\
  \midrule
  $\mathcal{Q}^-_r \Lambda^0$ & \texttt{Lagrange} & \texttt{Lagrange}  \\
   $\mathcal{Q}^-_r \Lambda^1$ & \texttt{RTCE} or \texttt{RTCF} & \texttt{NCE} \\
   $\mathcal{Q}^-_r \Lambda^2$ & \texttt{DQ} & \texttt{NCF}  \\
  $\mathcal{Q}^-_r \Lambda^3$ & \texttt{-} & \texttt{DQ}\\
  \midrule
  $\mathcal{S}^-_r \Lambda^0 $ & \texttt{S} & \texttt{S}  \\
  $\mathcal{S}^-_r \Lambda^1$ & \texttt{SminusCurl} or \texttt{SminusDiv} & \texttt{SminusCurl}  \\
  $\mathcal{S}^-_r \Lambda^2$ & \texttt{DPC} & \texttt{SminusDiv}\\
  $\mathcal{S}^-_r \Lambda^3$ & \texttt{-} & \texttt{DPC}\\
  \bottomrule
\end{tabular}
\end{table}

The newly supported elements, mapping FEEC spaces onto names in UFL are shown in the lower half of \cref{tab:FiredrakeNames}. Modifying the code from \cref{lst:mixedPoissonCode} to utilize trimmed serendipity spaces rather than tensor product spaces is then simply a case of replacing the element names in the \texttt{FunctionSpace} definitions with appropriate trimmed space names. Concretely, the new \texttt{FunctionSpace} definitions are shown in \cref{lst:pde_using_trimmed_serendipity}, the rest of the code remains unchanged.
\begin{lstlisting}[float=htbp,caption={Setting up Firedrake to use the trimmed serendipity elements in a mixed Poisson problem in 3D.}, label={lst:pde_using_trimmed_serendipity}, numbers=left, firstnumber=3, xleftmargin=20pt,  xrightmargin=20pt]
...
hDivSpace = FunctionSpace(mesh, "SminusDiv", polyDegree)
l2Space = FunctionSpace(mesh, "DPC", polyDegree - 1)
...
\end{lstlisting}

  \section{Experiments}
    
The following experiments show the benefits and costs of using trimmed serendipity elements in comparison to traditional tensor product elements.  We first present a basic projection example as a means of confirming approximation properties of our elements.  Next, we present results on a primal Poisson problem (to test $H^1$ elements), a mixed Poisson problem (to test \hdiv and $L^2$ elements), and a cavity resonator problem (to test \hcurl elements).  Since the \hcurl and \hdiv elements are only a rotation of the DOFs on the edges of elements in 2D, we do not give an experiment using \hcurl elements in 2D. 

The experiments were performed either on a single cluster compute node with 2x AMD EPYC 7642 48-core (Rome) processors (2.4GHz) and 512GB of memory running CentOS 7 or a similar node with 3TB of memory.  The 2D experiments were all done using the 512GB node, while in 3D, the 4th order experiments were run on the 3TB node.  Each job was run by submitting a SLURM script that requested one node in isolation to ensure no other jobs were running at the same time.  We utilized on-node parallelism by requesting a full node and executing the jobs with \texttt{mpirun} set to use 24 processes, with a few exceptions for smaller cases that will be pointed out later.  Timing data was collected after first performing a dry run of the code to warm the cache and then taking the minimum time over three subsequent runs.  The timing results presented here depend upon the solver choice, and changing that may give different results illustrating the relative efficiency between $\mathcal{Q}^-$ and $\mathcal{S}^-$.  

For simplicity, our numerical experiments all use a sparse direct solver.  We expect the multigrid theory in \citet{AFW2000,AFW2006} to carry over from existing $\mathcal{Q}^-$ spaces to $\mathcal{S}^-$ spaces.  Optimal smoothers require aggregating degrees of freedom for vertex patches, and we anticipate that the reduction in local dimensionality that trimmed serendipity spaces offer will be beneficial in these contexts as well.

\subsection{Projection}
  
We solve an $L^2$ projection problem to give a baseline accuracy test for the elements. Given either the unit square or unit cube as our domain of integration on definite integrals, we compute the projection of $f$ into the function space $V$ by using a discretization of the problem: find $u\in V \subset$ \hcurl such that
\begin{equation*}
  \int_\Omega u \cdot v \text{ d}x = \int_\Omega g \cdot v \text{ d}x, \text{ where } g = \nabla f
\end{equation*}
for all $v \in V$.  For our experiments, we choose $V$ to be \hcurl spaces for the proper dimension (one of the spaces \texttt{RTCE}, \texttt{NCE} or \texttt{SminusCurl}) and $f = \text{sin}(\pi x)\text{sin}(\pi y)$ or $f=\text{sin}(\pi x)\text{sin}(\pi y)\text{sin}(\pi z)$ for two or three dimensions respectively.
Recall that the trimmed serendipity elements are denoted with $\mathcal{S}_r^-$ and the tensor product elements with $\mathcal{Q}^-_r$.  In Firedrake, the trimmed serendipity elements use the label \texttt{SminusCurl} or \texttt{SminusDiv} for 1-forms in 2D, depending upon the orientation of the edge DOFs, and in 3D, they represent the 1- and 2-forms respectively.  The tensor product elements use the labels \texttt{RTCE} (or \texttt{RTCF}) and  \texttt{NCE} for 1-forms in 2D and 3D, while the 2-forms in 3D are  \texttt{NCF}.  To solve the projection problem, we use a Galerkin $L^2$ projection into $V$.

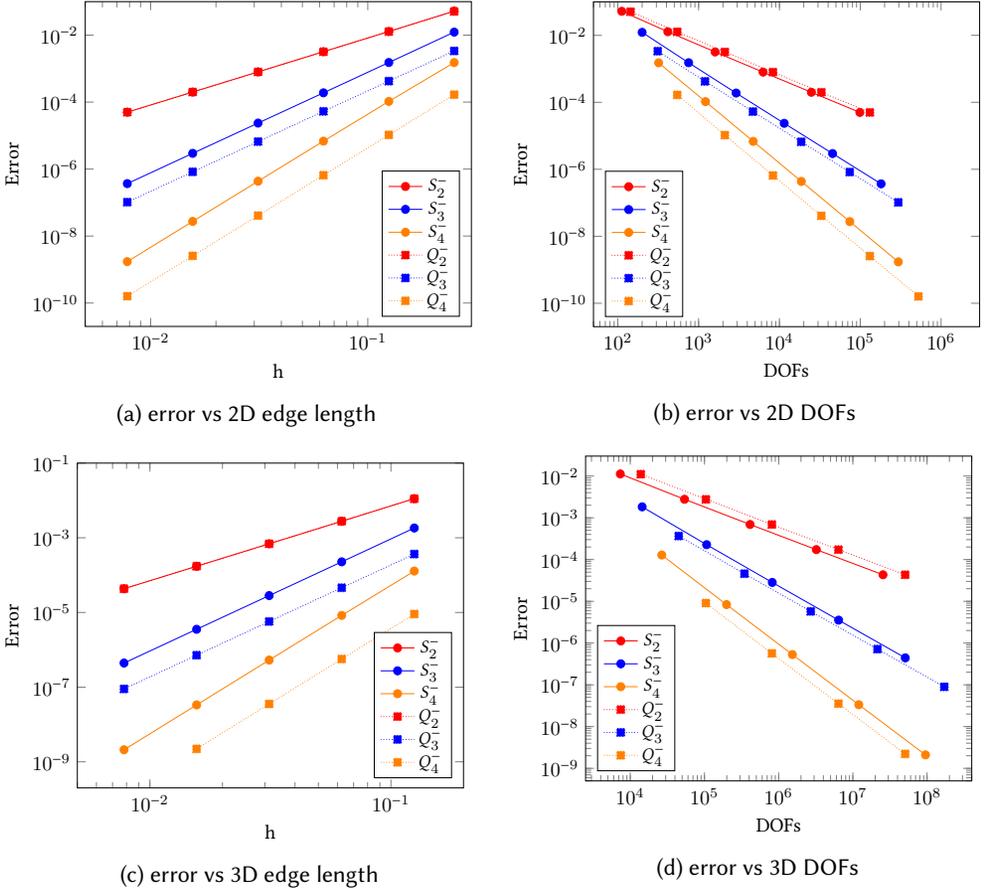
\begin{figure}[htbp]
  \centering
  \begin{subfigure}[h]{0.48\textwidth}
    \begin{tikzpicture}[scale=0.75]
      \begin{loglogaxis}[xlabel={h}, ylabel={Error},
             ylabel near ticks, ymax=1e-1, ymin=2e-11, xmax=0.3, xmin=0.5e-2,
             legend pos=south east, legend style={font=\small} ,
             cycle list name=mycolor]
        \addplot
        table [x=h,y=Error, col sep=comma]{csvs/ProjectionSminus2dO2.csv};
        \addlegendentry{$S^-_2$}
        \addplot 
        table [x=h,y=Error, col sep=comma]{csvs/ProjectionSminus2dO3.csv};
        \addlegendentry{$S^-_3$ }
        \addplot
        table [x=h,y=Error, col sep=comma]{csvs/ProjectionSminus2dO4.csv};
        \addlegendentry{$S^-_4$}
        \addplot
        table [x=h,y=Error, col sep=comma]{csvs/ProjectionRTCE2dO2.csv};
        \addlegendentry{$Q^-_2$}
        \addplot
        table [x=h,y=Error, col sep=comma]{csvs/ProjectionRTCE2dO3.csv};
        \addlegendentry{$Q^-_3$}
        \addplot[densely dotted, orange, mark=square*]
        table [x=h,y=Error, col sep=comma]{csvs/ProjectionRTCE2dO4.csv};
        \addlegendentry{$Q^-_4$}
      \end{loglogaxis}
      \end{tikzpicture}
      \caption{error vs 2D edge length \label{fig:2dProjectionH}}
    \end{subfigure}
    \begin{subfigure}[h]{0.48\textwidth}
      \begin{tikzpicture}[scale=0.75]
      \begin{loglogaxis}[xlabel={DOFs}, ylabel={Error},
             ylabel near ticks, ymax=1e-1, ymin=2e-11, xmax=3e6, xmin=0.5e2,
             legend pos=south west, legend style={font=\small} ,
             cycle list name=mycolor]
        \addplot
        table [x=Dofs,y=Error, col sep=comma]{csvs/ProjectionSminus2dO2.csv};
        \addlegendentry{$S^-_2$}
        \addplot
        table [x=Dofs,y=Error, col sep=comma]{csvs/ProjectionSminus2dO3.csv};
        \addlegendentry{$S^-_3$ }
        \addplot
        table [x=Dofs,y=Error, col sep=comma]{csvs/ProjectionSminus2dO4.csv};
        \addlegendentry{$S^-_4$}
        \addplot
        table [x=Dofs,y=Error, col sep=comma]{csvs/ProjectionRTCE2dO2.csv};
        \addlegendentry{$Q^-_2$}
        \addplot
        table [x=Dofs,y=Error, col sep=comma]{csvs/ProjectionRTCE2dO3.csv};
        \addlegendentry{$Q^-_3$}
        \addplot[densely dotted, orange, mark=square*]
        table [x=Dofs,y=Error, col sep=comma]{csvs/ProjectionRTCE2dO4.csv};
        \addlegendentry{$Q^-_4$}
      \end{loglogaxis}
      \end{tikzpicture}
      \caption{error vs 2D DOFs \label{fig:2dProjection}}
  \end{subfigure} \\[0.5\baselineskip]
  \begin{subfigure}[h]{0.48\textwidth}
    \begin{tikzpicture}[scale=0.75]
      \begin{loglogaxis}[xlabel={h}, ylabel={Error},
             ylabel near ticks, ymax=1e-1, ymin=2e-10, xmax=0.2, xmin=0.5e-2,
             legend pos=south east, legend style={font=\small} ,
             cycle list name=mycolor]
        \addplot
        table [x=h,y=Error, col sep=comma]{csvs/ProjectionSminus3dO2.csv};
        \addlegendentry{$S^-_2$}
        \addplot
        table [x=h,y=Error, col sep=comma]{csvs/ProjectionSminus3dO3.csv};
        \addlegendentry{$S^-_3$ }
        \addplot
        table [x=h,y=Error, col sep=comma]{csvs/ProjectionSminus3dO4.csv};
        \addlegendentry{$S^-_4$}
        \addplot
        table [x=h,y=Error, col sep=comma]{csvs/ProjectionNCE3dO2.csv};
        \addlegendentry{$Q^-_2$}
        \addplot
        table [x=h,y=Error, col sep=comma]{csvs/ProjectionNCE3dO3.csv};
        \addlegendentry{$Q^-_3$}
        \addplot[densely dotted, orange, mark=square*]
        table [x=h,y=Error, col sep=comma]{csvs/ProjectionNCE3dO4.csv};
        \addlegendentry{$Q^-_4$}
      \end{loglogaxis}
      \end{tikzpicture}
    \caption{error vs 3D edge length\label{fig:3dProjectionH}}
  \end{subfigure}
  \begin{subfigure}[h]{0.48\textwidth}
  \begin{tikzpicture}[scale=0.75]
      \begin{loglogaxis}[xlabel={DOFs}, ylabel={Error},
             ylabel near ticks, ymax=3e-2, ymin=0.5e-9, xmax=4e8, xmin=0.25e4,
             legend pos=south west, legend style={font=\small} ,
             cycle list name=mycolor]
        \addplot
        table [x=Dofs,y=Error, col sep=comma]{csvs/ProjectionSminus3dO2.csv};
        \addlegendentry{$S^-_2$}
        \addplot
        table [x=Dofs,y=Error, col sep=comma]{csvs/ProjectionSminus3dO3.csv};
        \addlegendentry{$S^-_3$ }
        \addplot
        table [x=Dofs,y=Error, col sep=comma]{csvs/ProjectionSminus3dO4.csv};
        \addlegendentry{$S^-_4$}
        \addplot
        table [x=Dofs,y=Error, col sep=comma]{csvs/ProjectionNCE3dO2.csv};
        \addlegendentry{$Q^-_2$}
        \addplot
        table [x=Dofs,y=Error, col sep=comma]{csvs/ProjectionNCE3dO3.csv};
        \addlegendentry{$Q^-_3$}
        \addplot[densely dotted, orange, mark=square*]
        table [x=Dofs,y=Error, col sep=comma]{csvs/ProjectionNCE3dO4.csv};
        \addlegendentry{$Q^-_4$}
      \end{loglogaxis}
      \end{tikzpicture}
    \caption{error vs 3D DOFs\label{fig:3dProjectionDofs}}
  \end{subfigure}
  \caption{Results of solving an $L^2$ projection problem using trimmed serendipity and tensor product \hcurl elements in both 2D (top row) and 3D (bottom row).  Experiments were ran for $k$-forms for $k=0,1,2,3$ in 2D and 3D, however only the $1$-forms in 3D are displayed here.  In every case, the trimmed serendipity element trendline has the same slope as its tensor product counterpart, as expected by the theory. In $(a)$ and $(c)$,  the order $2$ elements only show one trendline as they differ only slightly in error values and have the same edge lengths. \label{Projections}}
  
  
  \Description{Experimental convergence analysis on a projection problem using tensor product and trimmed serendipity 1-form elements.}
\end{figure}

The goal of the projection problem is to establish that the elements attain the mathematically expected $L^2$ rates of approximation as mesh size decreases, as well as comparing relative efficiencies of $\mathcal{S}^-$ and $\mathcal{Q}^-$.  For this, we create a mesh on $[0,1]^n$ of uniformly sized squares or cubes, where we refine the mesh from $N=4$ squares (or cubes) in each row and column to $N=128$ squares in each row and column (or $N=64$ cubes).  This results in a mesh with $N^2$ or $N^3$ squares or cubes respectively.  For the following results, we will use $h = \frac{1}{N}$ since the mesh elements are all uniformly sized. 
We employ both trimmed serendipity elements and tensor product elements on each mesh and record the $L^2$ error in each case.  In \cref{Projections}, we report the $L^2$ error in terms of the classical measure of maximum edge length ($h$) as well as total number of  DOFs.

The expectation is that $\mathcal{S}^-_r$ and $\mathcal{Q}^-_r$ converge at the same rate with respect to $h$, which is confirmed in the plots by parallel trendlines.  These parallel trendlines can be seen in each of the projection plots.  The overall results from the projection problem show that tensor product and trimmed serendipity elements give similar levels of error.  

For the $r=2$ case, the trimmed serendipity elements achieve a better accuracy while requiring fewer DOFs.  Therefore in this low order case, trimmed serendipity elements would be beneficial to use.  In the case of $r=3,4$ the trendlines for $\mathcal{S}^-_r$ are above the trendlines for $\mathcal{Q}^-_r$.  However, considering the elements $\mathcal{Q}^-_3$ and $\mathcal{S}^-_4$, we see that $\mathcal{Q}^-_3$ requires approximately $1.71 \times 10^8$ DOFs and $\mathcal{S}^-_4$ at the same mesh refinement requires $0.95 \times 10^8$ DOFs.  While $\mathcal{Q}^-_3$ attains an error of $8.96 \times 10^{-8}$, the $\mathcal{S}^-_4$ elements attain an error of $2.1 \times 10^{-9}$.

\subsection{The Poisson Problem}
In this section we discuss results for both the primal formulation and the mixed formulation of the Poisson problem.  We solve the primal weak formulation described below on a unit square domain $\Omega$ for $u \in U$:
\begin{align*}
    -\Delta u &= f \\
    \displaystyle u\vert_{\partial \Omega} &= 0
\end{align*}
where $f(x,y) = 2\pi^2\text{sin}(\pi x)\text{sin}(\pi y) $, yielding the solution $u(x,y) = \text{sin}(\pi x)\text{sin}(\pi y)$. In 3D, we can extend this to $f(x,y,z) = 3\pi^2\text{sin}(\pi x)\text{sin}(\pi y)\text{sin}(\pi z)$ with the solution $u(x,y,z) = \text{sin}(\pi x)\text{sin}(\pi y)\text{sin}(\pi z)$ on the unit cube.  The primal weak formulation of the Poisson equation is as follows: find $u \in V:= H^1(\Omega)$ such that:
\begin{equation*}
    \int_\Omega \nabla u \cdot \nabla v \text{ d}x = \int_\Omega f v \text{ d}x,\qquad \text{for all $v \in V$}.
\end{equation*}
Accordingly, the primal formulation employs $H^1$ elements, using \texttt{S} for $\mathcal{S}^-_r \Lambda^0$ and \texttt{Lagrange} for $\mathcal{Q}^-_r \Lambda^0$.

The mixed formulation of the Poisson problem introduces an intermediate variable, $\sigma$, which is solved for simultaneously.
Formally, this is: find $\sigma\in$~\hdiv and $u\in L^2$ such that:
\begin{align}
     \sigma - \nabla u &= 0 \notag \\
     \nabla \cdot \sigma &= -f \notag \\
     u\vert_{\partial \Omega} &= 0 \notag
\end{align}
In a similar fashion as for the primal formulation, we can create the mixed formulation of the Poisson problem that we saw in \cref{eq:MixedPoisson}.

These equations are discretized and solved using a suitable \textit{pair} of finite elements -- one of \hdiv type and one of $L^2$ type.  We use $(\calQ_r^-\Lambda^{n-1},\calQ_r^-\Lambda^n)$ and $(\calS_r^-\Lambda^{n-1},\calS_r^-\Lambda^n)$ for dimensions $n=2$ and $n=3$.  Note that the mathematical notation here calls for $\mathcal{Q}^-_r \Lambda^{n-1}$ to be paired with $\mathcal{Q}^-_r \Lambda^n$, but the code notation will require setting the degree of the $L^2$ element one below the degree of the \hdiv element, and similar with the trimmed serendipity elements.
For both the primal Poisson and mixed Poisson problems, we solve the system using MUMPS \citep{MUMPS01,MUMPS02} with a sparse direct LU factorization using iterative refinement in order to attain high accuracy in the solvers and allow us to focus on analyzing the elements instead of confounding variables.  At high degree and fine mesh resolutions, we noticed that both the tensor product and trimmed serendipity elements would hit a floor in error values that was unexpected. At lower degrees or coarser meshes, this was unnecessary, but we chose to keep the solver parameters the same to be consistent throughout the experiments. The exact details of the MUMPS configuration can be found both in the zenodo link and in the appendix in order for results to be reproducible.


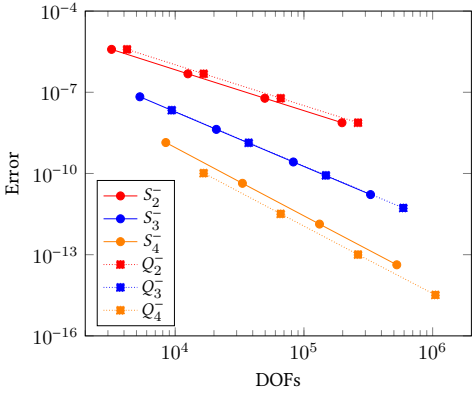
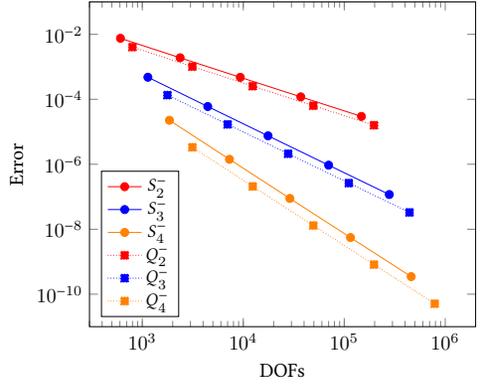
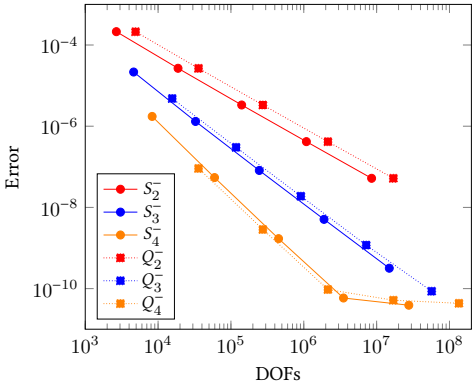
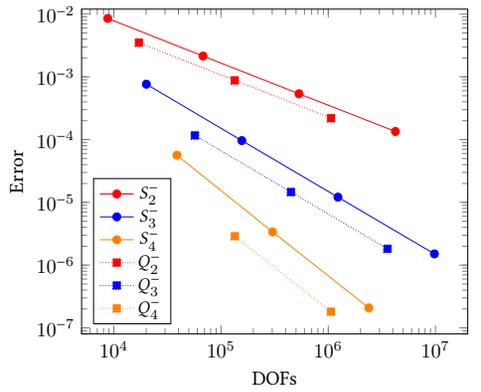
\begin{figure}[htbp]
  \centering
  \begin{subfigure}[h]{0.48\textwidth}
      \begin{tikzpicture}[scale=0.75]
      \begin{loglogaxis}[xlabel={DOFs}, ylabel={Error},
             ylabel near ticks, ymax=1e-4, ymin=1e-16, xmax=2e6, xmin=2e3,
             legend pos=south west, legend style={font=\small} ,
             cycle list name=mycolor]
        \addplot
        table [x=Dofs,y=Error, col sep=comma]{csvs/PrimalPoissonSerendipity2dO2.csv};
        \addlegendentry{$S^-_2$}
        \addplot
        table [x=Dofs,y=Error, col sep=comma]{csvs/PrimalPoissonSerendipity2dO3.csv};
        \addlegendentry{$S^-_3$ }
        \addplot
        table [x=Dofs,y=Error, col sep=comma]{csvs/PrimalPoissonSerendipity2dO4.csv};
        \addlegendentry{$S^-_4$}
        \addplot
        table [x=Dofs,y=Error, col sep=comma]{csvs/PrimalPoissonLagrange2dO2.csv};
        \addlegendentry{$Q^-_2$}
        \addplot
        table [x=Dofs,y=Error, col sep=comma]{csvs/PrimalPoissonLagrange2dO3.csv};
        \addlegendentry{$Q^-_3$}
        \addplot[densely dotted, orange, mark=square*]
        table [x=Dofs,y=Error, col sep=comma]{csvs/PrimalPoissonLagrange2dO4.csv};
        \addlegendentry{$Q^-_4$}
      \end{loglogaxis}
      \end{tikzpicture}
      \caption{2D primal Poisson convergence analysis. \label{fig:2dPrimalDofs}}
  \end{subfigure}
  \begin{subfigure}[h]{0.48\textwidth}
    \begin{tikzpicture}[scale=0.75]
      \begin{loglogaxis}[xlabel={DOFs}, ylabel={Error},
             ylabel near ticks, ymax=1e-1, ymin=1e-11, xmax=2e6, xmin=3e2,
             legend pos=south west, legend style={font=\small} ,
             cycle list name=mycolor]
        \addplot
        table [x=Dofs,y=Error, col sep=comma]{csvs/MixedPoissonSminus2dO2.csv};
        \addlegendentry{$S^-_2$}
        \addplot
        table [x=Dofs,y=Error, col sep=comma]{csvs/MixedPoissonSminus2dO3.csv};
        \addlegendentry{$S^-_3$ }
        \addplot
        table [x=Dofs,y=Error, col sep=comma]{csvs/MixedPoissonSminus2dO4.csv};
        \addlegendentry{$S^-_4$}
        \addplot
        table [x=Dofs,y=Error, col sep=comma]{csvs/MixedPoissonRTCF2dO2.csv};
        \addlegendentry{$Q^-_2$}
        \addplot
        table [x=Dofs,y=Error, col sep=comma]{csvs/MixedPoissonRTCF2dO3.csv};
        \addlegendentry{$Q^-_3$}
        \addplot[densely dotted, orange, mark=square*]
        table [x=Dofs,y=Error, col sep=comma]{csvs/MixedPoissonRTCF2dO4.csv};
        \addlegendentry{$Q^-_4$}
      \end{loglogaxis}
      \end{tikzpicture}
    \caption{2D mixed Poisson convergence analysis. \label{fig:2dMixedDofsError}}
  \end{subfigure} \\[0.5\baselineskip]
  \begin{subfigure}[h]{0.48\textwidth}
    \begin{tikzpicture}[scale=0.75]
      \begin{loglogaxis}[xlabel={DOFs}, ylabel={Error},
             ylabel near ticks, ymax=1e-3, ymin=1e-11, xmax=2e8, xmin=1e3,
             legend pos=south west, legend style={font=\small} ,
             cycle list name=mycolor]
        \addplot
        table [x=Dofs,y=Error, col sep=comma]{csvs/PrimalPoissonSerendipity3dO2.csv};
        \addlegendentry{$S^-_2$}
        \addplot
        table [x=Dofs,y=Error, col sep=comma]{csvs/PrimalPoissonSerendipity3dO3.csv};
        \addlegendentry{$S^-_3$ }
        \addplot
        table [x=Dofs,y=Error, col sep=comma]{csvs/PrimalPoissonSerendipity3dO4.csv};
        \addlegendentry{$S^-_4$}
        \addplot
        table [x=Dofs,y=Error, col sep=comma]{csvs/PrimalPoissonLagrange3dO2.csv};
        \addlegendentry{$Q^-_2$}
        \addplot
        table [x=Dofs,y=Error, col sep=comma]{csvs/PrimalPoissonLagrange3dO3.csv};
        \addlegendentry{$Q^-_3$}
        \addplot[densely dotted, orange, mark=square*]
        table [x=Dofs,y=Error, col sep=comma]{csvs/PrimalPoissonLagrange3dO4.csv};
        \addlegendentry{$Q^-_4$}
      \end{loglogaxis}
      \end{tikzpicture}
    \caption{3D primal Poisson convergence analysis. \label{fig:3dPrimalDofsError}}
  \end{subfigure}
  \begin{subfigure}[h]{0.48\textwidth}
    \begin{tikzpicture}[scale=0.75]
      \begin{loglogaxis}[xlabel={DOFs}, ylabel={Error},
             ylabel near ticks, ymax=1.2e-2, ymin=0.8e-7, xmax=2e7, xmin=5e3,
             legend pos=south west, legend style={font=\small} ,
             cycle list name=mycolor]
        \addplot
        table [x=Dofs,y=Error, col sep=comma]{csvs/MixedPoissonSminus3dO2.csv};
        \addlegendentry{$S^-_2$}
        \addplot
        table [x=Dofs,y=Error, col sep=comma]{csvs/MixedPoissonSminus3dO3.csv};
        \addlegendentry{$S^-_3$ }
        \addplot
        table [x=Dofs,y=Error, col sep=comma]{csvs/MixedPoissonSminus3dO4.csv};
        \addlegendentry{$S^-_4$}
        \addplot
        table [x=Dofs,y=Error, col sep=comma]{csvs/MixedPoissonNCF3dO2.csv};
        \addlegendentry{$Q^-_2$}
        \addplot
        table [x=Dofs,y=Error, col sep=comma]{csvs/MixedPoissonNCF3dO3.csv};
        \addlegendentry{$Q^-_3$}
        \addplot[densely dotted, orange, mark=square*]
        table [x=Dofs,y=Error, col sep=comma]{csvs/MixedPoissonNCF3dO4.csv};
        \addlegendentry{$Q^-_4$}
      \end{loglogaxis}
      \end{tikzpicture}
      \caption{3D mixed Poisson convergence analysis. \label{fig:3dMixedDofsError}}
  \end{subfigure}
  \caption{A convergence analysis of the primal and mixed Poisson problems in 2D and 3D.  Error here is calculated as the $L^2$ error between the exact solution and the approximate using the corresponding finite element space.}
\label{fig:PrimalMixedErrorAnalysis}
\Description{Experimental convergence analysis for the primal and mixed poisson problems using trimmed serendipity and tensor product element.  Uses $H^1$ elements for the primal problem and \hdiv and $L^2$ elements for the mixed problem.  }
\end{figure}

The empirical convergence results for the primal and mixed formulations of the Poisson problems can be seen in \cref{fig:2dPrimalDofs,fig:2dMixedDofsError,fig:3dPrimalDofsError,fig:3dMixedDofsError}.  In each of the subfigures of \cref{fig:PrimalMixedErrorAnalysis}, we see that independent of the performance of each element, the $\mathcal{S}^-_r$ and $\mathcal{Q}^-_r$ have parallel trendlines, indicating that they have the same overall convergence rate.  For the primal formulation in 2D and 3D, the trimmed serendipity elements perform similar to the tensor product elements for orders $r=2,3$.  Furthermore, comparing the elements via DOFs as in the projection problem yields another instance where we see that using $\mathcal{S}^-_{3}$ instead of $\mathcal{Q}^-_2$ will attain a higher accuracy for essentially the same number of DOFs.

\begin{figure}[htbp]
  \centering
  \begin{subfigure}[h]{0.48\textwidth}
    \begin{tikzpicture}[scale=0.75]
      \begin{loglogaxis}[xlabel={DOFs}, ylabel={Time},
             ylabel near ticks, ymax=1.e+4, ymin=0.5e-1, xmax=2.e8, xmin=1.e+3,
             legend pos=north west, legend style={font=\small} ,
             cycle list name=mycolor]
        \addplot
        table [x=Dofs,y=Time, col sep=comma]{csvs/PrimalPoissonSerendipity3dO2.csv};
        \addlegendentry{$S^-_2$}
        \addplot
        table [x=Dofs,y=Time, col sep=comma]{csvs/PrimalPoissonSerendipity3dO3.csv};
        \addlegendentry{$S^-_3$ }
        \addplot
        table [x=Dofs,y=Time, col sep=comma]{csvs/PrimalPoissonSerendipity3dO4.csv};
        \addlegendentry{$S^-_4$}
        \addplot
        table [x=Dofs,y=Time, col sep=comma]{csvs/PrimalPoissonLagrange3dO2.csv};
        \addlegendentry{$Q^-_2$}
        \addplot
        table [x=Dofs,y=Time, col sep=comma]{csvs/PrimalPoissonLagrange3dO3.csv};
        \addlegendentry{$Q^-_3$}
        \addplot[densely dotted, orange, mark=square*]
        table [x=Dofs,y=Time, col sep=comma]{csvs/PrimalPoissonLagrange3dO4.csv};
        \addlegendentry{$Q^-_4$}
      \end{loglogaxis}
      \end{tikzpicture}
      \caption{3D primal Poisson Time vs DOFs \label{fig:3dPrimalDofsTime}}
  \end{subfigure}
  \begin{subfigure}[h]{0.48\textwidth}
    \begin{tikzpicture}[scale=0.75]
      \begin{loglogaxis}[xlabel={Time}, ylabel={Error},
             ylabel near ticks, ymax=1e-3, ymin=1e-11, xmax=1e4, xmin=0.5e-1,
             legend pos=north east, legend style={font=\small} ,
             cycle list name=mycolor]
        \addplot
        table [x=Time,y=Error, col sep=comma]{csvs/PrimalPoissonSerendipity3dO2.csv};
        \addlegendentry{$S^-_2$}
        \addplot
        table [x=Time,y=Error, col sep=comma]{csvs/PrimalPoissonSerendipity3dO3.csv};
        \addlegendentry{$S^-_3$ }
        \addplot
        table [x=Time,y=Error, col sep=comma]{csvs/PrimalPoissonSerendipity3dO4.csv};
        \addlegendentry{$S^-_4$}
        \addplot
        table [x=Time,y=Error, col sep=comma]{csvs/PrimalPoissonLagrange3dO2.csv};
        \addlegendentry{$Q^-_2$}
        \addplot
        table [x=Time,y=Error, col sep=comma]{csvs/PrimalPoissonLagrange3dO3.csv};
        \addlegendentry{$Q^-_3$}
        \addplot[densely dotted, orange, mark=square*]
        table [x=Time,y=Error, col sep=comma]{csvs/PrimalPoissonLagrange3dO4.csv};
        \addlegendentry{$Q^-_4$}
      \end{loglogaxis}
      \end{tikzpicture}
      \caption{3D primal Poisson Error vs Time \label{fig:3dPrimalTimeError}}
  \end{subfigure}\\[0.5\baselineskip]
  \begin{subfigure}[h]{0.48\textwidth}
    \begin{tikzpicture}[scale=0.75]
      \begin{loglogaxis}[xlabel={DOFs}, ylabel={Time},
             ylabel near ticks, ymax=2e4, ymin=1e-1, xmax=2e7, xmin=0.5e4,
             legend pos=north west, legend style={font=\small} ,
             cycle list name=mycolor]
        \addplot
        table [x=Dofs,y=Time, col sep=comma]{csvs/MixedPoissonSminus3dO2.csv};
        \addlegendentry{$S^-_2$}
        \addplot
        table [x=Dofs,y=Time, col sep=comma]{csvs/MixedPoissonSminus3dO3.csv};
        \addlegendentry{$S^-_3$ }
        \addplot
        table [x=Dofs,y=Time, col sep=comma]{csvs/MixedPoissonSminus3dO4.csv};
        \addlegendentry{$S^-_4$}
        \addplot
        table [x=Dofs,y=Time, col sep=comma]{csvs/MixedPoissonNCF3dO2.csv};
        \addlegendentry{$Q^-_2$}
        \addplot
        table [x=Dofs,y=Time, col sep=comma]{csvs/MixedPoissonNCF3dO3.csv};
        \addlegendentry{$Q^-_3$}
        \addplot[densely dotted, orange, mark=square*]
        table [x=Dofs,y=Time, col sep=comma]{csvs/MixedPoissonNCF3dO4.csv};
        \addlegendentry{$Q^-_4$}
      \end{loglogaxis}
      \end{tikzpicture}
    \caption{3D mixed Poisson Time vs DOFs \label{fig:3dMixedDofsTime}}
  \end{subfigure}
  \begin{subfigure}[h]{0.48\textwidth}
    \begin{tikzpicture}[scale=0.75]
      \begin{loglogaxis}[xlabel={Time}, ylabel={Error},
             ylabel near ticks, ymax=1e-1, ymin=5e-8, xmax=1e5, xmin=0.5e-1,
             legend pos=north east, legend style={font=\small} ,
             cycle list name=mycolor]
        \addplot
        table [x=Time,y=Error, col sep=comma]{csvs/MixedPoissonSminus3dO2.csv};
        \addlegendentry{$S^-_2$}
        \addplot
        table [x=Time,y=Error, col sep=comma]{csvs/MixedPoissonSminus3dO3.csv};
        \addlegendentry{$S^-_3$ }
        \addplot
        table [x=Time,y=Error, col sep=comma]{csvs/MixedPoissonSminus3dO4.csv};
        \addlegendentry{$S^-_4$}
        \addplot
        table [x=Time,y=Error, col sep=comma]{csvs/MixedPoissonNCF3dO2.csv};
        \addlegendentry{$Q^-_2$}
        \addplot
        table [x=Time,y=Error, col sep=comma]{csvs/MixedPoissonNCF3dO3.csv};
        \addlegendentry{$Q^-_3$}
        \addplot[densely dotted, orange, mark=square*]
        table [x=Time,y=Error, col sep=comma]{csvs/MixedPoissonNCF3dO4.csv};
        \addlegendentry{$Q^-_4$}
      \end{loglogaxis}
      \end{tikzpicture}
    \caption{3D Mixed Poisson Error vs Time \label{fig:3dMixedTimeError}}
  \end{subfigure}
  \caption{Analyzing timing data for primal and mixed Poisson problems using trimmed serendipity and tensor product elements.  Error here is calculated as the $L^2$ error between the exact solution and the approximate solution found using the corresponding finite element space.\label{fig:PrimalMixedTimeAnalysis}}
\Description{Some timing results for the primal and mixed poisson problems using trimmed serendipity and tensor product element.  Uses $H^1$ elements for the primal problem and \hdiv and $L^2$ elements for the mixed problem.  }
\end{figure}
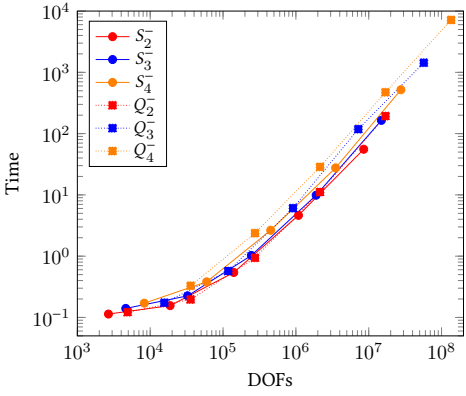
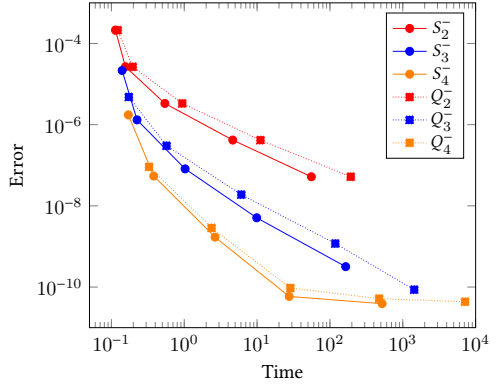
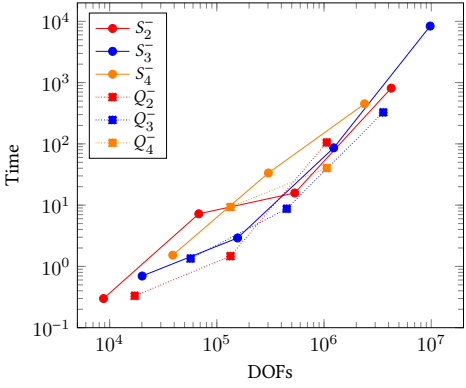
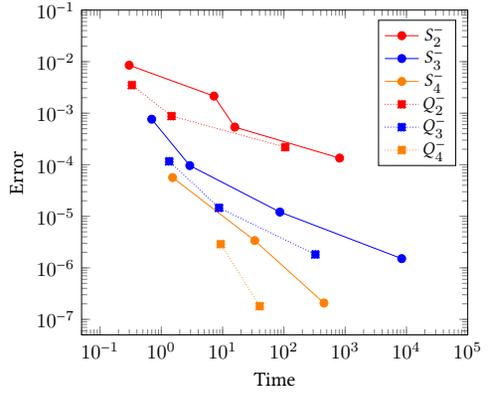

In \cref{fig:PrimalMixedTimeAnalysis} we analyze the timing data for computing the solutions to the primal and mixed formulations using trimmed serendipity and tensor product elements.  
As in the error vs DOFs graphs, we see good evidence in \cref{fig:3dPrimalDofsTime,fig:3dPrimalTimeError} that the trimmed serendipity and tensor product elements compute solutions at a similar speed based on the number of DOFs.
In \cref{fig:3dPrimalTimeError} we see that for a given error level, trimmed serendipity elements require less time.  The overall time required being dependent upon on the number of DOFs rather than the element type is seen again in  \cref{fig:3dMixedDofsTime}.  Further evidence of this is seen in \cref{fig:3dMixedTimeError}.  Similar to the previous analysis of DOFs vs error for the mixed formulation, the timing data here illustrates that attaining the extra accuracy from using $\mathcal{S}^-_3$ instead of $\mathcal{Q}^-_2$ does not invoke a larger time requirement. The sparsity of the matrices involved in the order 4 elements can be seen in \cref{tab:NNZ}.

\begin{table}[htbp]
  \centering
  \caption{Expressing the number of nonzero entries in the matrices used to compute solutions to the primal and mixed formulations of the Poisson problem.  The data shown here represents order 4 elements, where the meshes are either $128^2$ or $64^3$ depending on the dimension of the space.  The first row of each half indicates the number of nonzero entries, while the second row of each half indicates the proportion of the number of nonzero entries. \label{tab:NNZ}}
  \begin{tabular}{ c c c c }
    \toprule
\multicolumn{4}{c}{$\mathcal{Q}^-_4$ Elements} \\[0.3em]
 Primal $n=2$ & Primal $n=3$ & Mixed $n=2$ & Mixed $n=3$ \\ 
\midrule
381,825 & 143,992,308 & 2,096,704 & 989,178,624 \\  
\num{5.51e-06} & \num{4.99e-07} & \num{3.38e-06} & \num{2.18e-07}\\
\bottomrule
\end{tabular} \\[\baselineskip]
\begin{tabular}{ c c c c }
  \toprule
\multicolumn{4}{c}{$\mathcal{S}^-_4$ Elements} \\[0.3em]
 Primal $n=2$ & Primal $n=3$ & Mixed $n=2$ & Mixed $n=3$ \\ 
\midrule
156,625 & 17,148,900 & 848,624 & 107,771,596 \\  
\num{8.97e-06} & \num{1.39e-06} & \num{4.01e-06} & \num{2.98e-07}\\
\bottomrule
\end{tabular}
\end{table}

\subsection{Cavity Resonator}

The last numerical experiment that we give here is the cavity resonator problem, making use of the \hcurl elements in 3D.  We a pose a Maxwell eigenvalue problem on the domain $\Omega = [0,1]^3$ with perfectly conducting boundary conditions, yielding an eigenvalue problem where $\lambda$ represents a quantity proportional to the frequency squared of the time-harmonic electric field (i.e.\ eigenvalues) and $E$ represents the electric field (i.e.\ eigenfunctions):
\begin{align*}
    \nabla \cdot E &= 0 \text{ in } \Omega \\
    \nabla \times \nabla \times E &= \lambda E  \text{ in } \Omega \\
    E \times n &= 0 \text{ on } \partial \Omega.
\end{align*}
We consider the weak formulation of this problem (similar to \citet{fumio1987mixed}), where $\omega$ represents the resonances (i.e.\ eigenvalues) and $E$ represents the electric field (i.e.\ eigenfunctions):
\begin{equation*}
    \int_\Omega \big(\nabla \times F\big) \cdot \big(\nabla \times E \big) \text{ d}x = \omega^2 \int_\Omega F \cdot E \text{ d}x \text{ for all } F \in H_0(\text{curl}).
\end{equation*}
The exact eigenvalues follow the formula
\[ \omega^2 = m_1^2 + m_2^2 + m_3^2 \]
where $m_i \in \mathbb{N} \cup {0}$ and no more than one of $m_1, m_2, m_3$ may be equal to $0$ at a time \cite{rognes2010efficient}.

\begin{table}[htbp]
  \centering
  \caption{A comparison of how $ \mathcal{Q}^-_2$ and $\mathcal{S}^-_2$ finite elements solve the Maxwell cavity resonator eigenvalue problem, $\langle \text{curl}(F), \text{curl}(E) \rangle = \omega^2 \langle F, E \rangle$. An eigenvalue found with the same error multiple times was condensed to a single row.  Numbers in parentheses next to the actual eigenvalue are the number of times we found an approximation of the actual eigenvalue.  The columns labeled $N=4, 8, 16, 32$ are giving the approximate eigenvalues found on a mesh of size $N \times N \times N$.  The values in parentheses in these columns indicates the rate of convergence for that approximate eigenvalue.\label{tab:Eigenvalue}}
  \begin{tabular}{ c c c c c }
    \toprule
\multicolumn{5}{c}{$\mathcal{Q}^-_2$ \hcurl Elements} \\[0.3em]
Actual (Count) & N = 4 & N = 8 & N = 16 & N = 32 \\ 
\midrule
2 (3) &2.001024 & 2.000066 (3.96) & 2.000004 (4.04) & 2.0000003 (4.00) \\  
3 (2) & 3.001536 & 3.000098 (3.97) & 3.000006 (4.03) & 3.0000004 (4.02) \\
5 (4) & 5.030601 & 5.002081 (3.88)& 5.000133 (3.97) & 5.000008 (4.06) \\
6 (3) & 6.031114 & 6.002114 (3.88) & 6.000135 (3.97) &  6.000008 (4.08) \\
8 (0) & $-$ & $-$ & $-$ & $-$ \\
\midrule
DOF  & 1944 & 13872 & 104544 & 811200 \\
\midrule
EPS solve time per iteration & 0.01565225 & 0.0743845 & 1.0484236 & 7.6186526 \\
\bottomrule
\end{tabular} \\[\baselineskip]

\begin{tabular}{ c c c c c }
  \toprule
\multicolumn{5}{c}{$\mathcal{S}^-_2$ \hcurl Elements} \\[0.3em]
Actual (Count) & N = 4 & N = 8 & N = 16 & N = 32 \\ 
\midrule
2 (3) & 2.001092 & 2.000066 (4.05) & 2.000004 (4.04) & 2.000000 (4.00) \\  
3 (2) & 3.009018 & 3.000586 (3.94) & 3.000037 (3.99) & 3.000002 (4.21) \\
5 (3) & 5.032027 & 5.002097 (3.93)& 5.000133 (3.98) & 5.000008 (4.06) \\
5 (1) & 5.032027 & 5.002097 (3.93) & 5.000133 (3.98) & $-$ \\
6 (1) & 6.072012 & 6.004976 (3.86) & 6.000319 (3.96) & 6.000020 (4.00) \\
6 (1) & 6.072012 & 6.004976 (3.86) & 6.000319 (3.96) & 6.000024 (3.73)\\
6 (1) & $-$ & $-$ & 6.00038 & 6.000024 (3.98)\\
8 (1) & $-$ & $-$ & $-$ & 8.000017 \\
\midrule
DOF  & 1080 & 7344 & 53856 & 411840 \\
\midrule
EPS solve time per iteration & 0.01288725 & 0.0309768 & 0.401663 & 4.1996873 \\
\bottomrule
\end{tabular}
\end{table}
In \cref{tab:Eigenvalue}, we display the convergence rates of different eigenvalues when computing the eigenvalues with tensor product and trimmed serendipity elements in 3D.  The table is split into halves, the top half showing values from using $\mathcal{Q}^-$ \hcurl elements while the bottom half shows values from using the corresponding $\mathcal{S}^-$ elements.  Each half of the table has a row giving the DOFs in the mesh for each refinement level $N$ and a row giving the time per iteration that the solver required.  

Note that the convergence rates are computed by
\[r = \frac{\text{log}\bigg(\frac{\tilde{\lambda}_{i,N} - \lambda_{i,N}}{\tilde{\lambda}_{i,N+1} - \lambda_{i,N+1}} \bigg)}{\text{log}\bigg( \frac{h_N}{h_{N+1}} \bigg)}. \]
Based off earlier eigenvalue works \cite{boffi2010finite}, we expect the rate of convergence to be double the order of the finite element used to solve the problem.  This is reflected in the table well for both $\mathcal{S}^-$ and $\mathcal{Q}^-$ elements.  The ``$-$'' entries in the table indicate the eigenvalue solver did not find that specific eigenvalue in the allowed number of iterations; we set the solver to iterate a sufficient number of times to find the first 15 eigenvalue-eigenvector pairs.

The experiment was done by using Firedrake to create the mass and stiffness matrices as petsc4py objects \citep{petsc-user-ref,petsc-efficient,Dalcin2011}, then using slepc4py \citep{slepc-toms,slepc-manual} to do the eigenvalue analysis.  The eigenvalue analysis was done by computing an inverted shift to a target of $3.0$, then solving for $15$ eigenvalue-eigenvector pairs.  The SLEPc solve was done using the default (Krylov-Schur) solver with a tolerance level of $10^{-7}$.
We note that the eigensolver finds a varying number of spurious eigenvalues with value $1$.  These exist because Firedrake enforces strong boundary conditions by placing a $1$ on the diagonal and zeroing out the rows and columns, not due to the elements that we use or the SLEPc solver that is called.  We do not report these eigenvalues.

Since both elements attain the expected convergence rate, we focus on the rest of the results in the table.  Investigating the error in the eigenvalues in the chart compared to the exact values, we see that tensor product elements are able to get results that are up to an order of magnitude better near the target eigenvalue.  On the other hand, this loss of accuracy from using trimmed serendipity elements is offset by a reduction in required time to solve for the requested eigenvalues.  At every mesh refinement level, trimmed serendipity elements have nearly half the DOFs of tensor product elements, and correspondingly, require approximately half the time per iteration to solve for the eigenvalues (outside of the case $N=4$).  At higher orders, we expect that this will be even more exaggerated.  

Continuing the eigenvalue example, we used Firedrake and SLEPc to compute two eigenvalues, $\lambda = 3$ and $\lambda = 5$. We computed the eigenvalues at different orders of the elements, from $r=2$ to $r=5$, and kept the mesh constant at $16 \times 16 \times 16$.  The timing data was then collected by choosing the largest time required for any of the multiplicities of $3$ or $5$ that the solver found. 

\begin{figure}[htbp]
  \centering
  \begin{subfigure}[h]{0.48\textwidth}
  \begin{tikzpicture}[scale=0.75]
      \begin{loglogaxis}[xlabel={DOFs}, ylabel={Error},
             ylabel near ticks, ymax=1e-3, ymin=1e-14, xmax=3e6, xmin=3.5e4,
             legend pos=north east, legend style={font=\small}]
        \addplot[red, mark=*]
        table [x=Dofs,y=Error, col sep=comma]{csvs/SminusEigenvaluesValue3.csv};
        \addlegendentry{$S^- \lambda = 3$}
        \addplot[densely dotted, red, mark=square*]
        table [x=Dofs,y=Error, col sep=comma]{csvs/TensorEigenvaluesValue3.csv};
        \addlegendentry{$Q^- \lambda = 3$}
        \addplot[blue, mark=*] 
        table [x=Dofs,y=Error, col sep=comma]{csvs/SminusEigenvaluesValue5.csv};
        \addlegendentry{$S^- \lambda = 5$ }
        \addplot[densely dotted, blue, mark=square*]
        table [x=Dofs,y=Error, col sep=comma]{csvs/TensorEigenvaluesValue5.csv};
        \addlegendentry{$Q^- \lambda = 5$}
      \end{loglogaxis}
      \end{tikzpicture}
      \caption{Eigenvalue error analysis. \label{fig:EigenvalueDofsError}}
     \end{subfigure}
  \begin{subfigure}[h]{0.48\textwidth}
    \begin{tikzpicture}[scale=0.75]
      \begin{loglogaxis}[xlabel={Time}, ylabel={Error},
             ylabel near ticks, ymax=1e-3, ymin=1e-14, xmax=185, xmin=0.1,
             legend pos=north east, legend style={font=\small}]
        \addplot[red, mark=*]
        table [x=Time,y=Error, col sep=comma]{csvs/SminusEigenvaluesValue3.csv};
        \addlegendentry{$S^- \lambda = 3$}
        \addplot[densely dotted, red, mark=square*]
        table [x=Time,y=Error, col sep=comma]{csvs/TensorEigenvaluesValue3.csv};
        \addlegendentry{$Q^- \lambda = 3$}
        \addplot[blue, mark=*] 
        table [x=Time,y=Error, col sep=comma]{csvs/SminusEigenvaluesValue5.csv};
        \addlegendentry{$S^- \lambda = 5$ }
        \addplot[densely dotted, blue, mark=square*]
        table [x=Time,y=Error, col sep=comma]{csvs/TensorEigenvaluesValue5.csv};
        \addlegendentry{$Q^- \lambda = 5$}
      \end{loglogaxis}
      \end{tikzpicture}
    \caption{Eigenvalue time analysis.\label{fig:EigenvalueTimeError}}
  \end{subfigure}\\
  \caption{Results for solving for $\lambda = 3$ and $\lambda = 5$ using Firedrake and SLEPc by increasing the order from $2$ to $5$.  Error is calculated as the absolute value of the error between the actual eigenvalue and the approximated eigenvalue.\label{fig:EigenvalueAnalysis}}
\Description{Convergence and timing analysis for the Maxwell cavity eigenvalue problem using trimmed serendipity and tensor product \hcurl elements.}
\end{figure}
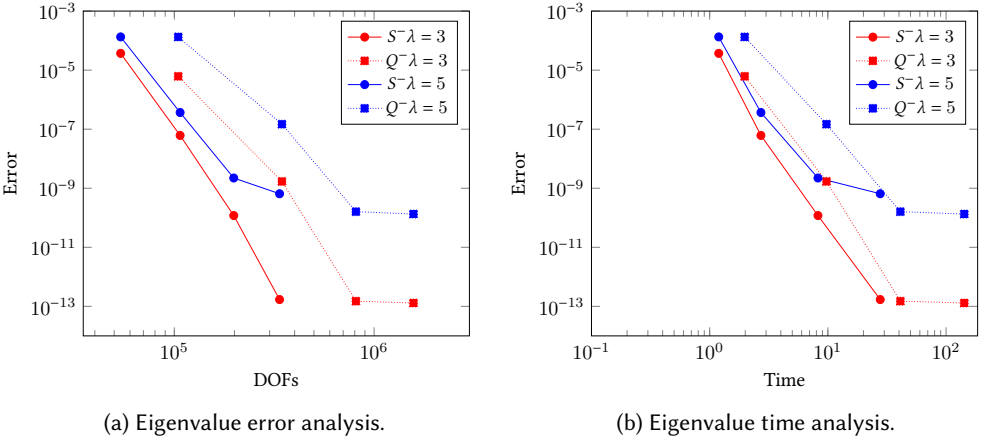

The error results shown in \cref{fig:EigenvalueDofsError} for the eigenvalue problem indicate that trimmed serendipity elements yield less error in the eigenvalues for the number DOFs required to compute them than the tensor product elements.  This is a change from the mixed formulation of the Poisson problem where the DOFs vs Error trendline for trimmed serendipity was generally above the trendline for tensor product elements.
The timing results in \cref{fig:EigenvalueTimeError} showed that the timing requirements for both trimmed serendipity and tensor product elements were similar, with trimmed serendipity generally requiring a little bit less time for a given error value.

\section{Discussion}

This implementation of trimmed serendipity elements gives a new method for computing the solution to a discretized PDE and has been tested on meshes of squares and cubes.  Completing the implementation of these elements within Firedrake by using the basis functions defined in \citet{gillette2019computational} is an illustration of Firedrake's modular capabilities for implementing new and unusual finite elements. 

The convergence studies done in each of the numerical experiments show that the trimmed serendipity elements can attain the theoretical rates of convergence that they were predicted to achieve.  While we only illustrate orders 2, 3, and 4 in 2D and 3D, our implementation of trimmed serendipity elements in Firedrake is designed to work in both 2D and 3D for arbitrary orders $r$.  

In comparison to tensor product elements $\mathcal{Q}^-_r$, we make a choice when using trimmed serendipity elements $\mathcal{S}^-_r$ to lower accuracy in return for less computation, both in terms of DOFs and time required.  At low orders the choice to use trimmed serendipity elements could actually reduce the error per DOF, as we saw in the primal formulation of the Poisson problem, where the trendlines for DOFs vs Error for trimmed serendipity elements were below the trendlines for the tensor product elements.  In the mixed formulation case however, the opposite was true, and the trendlines for the tensor product elements were below the trendlines for the trimmed serendipity elements.

Rather than comparing in terms of approximation order, it can also be beneficial to compare the two elements based off of the DOFs that they require.  Consider the 3D mixed formulation of the Poisson problem while focusing on DOFs vs Error given in \cref{fig:3dMixedDofsError}.  The tensor product elements $\mathcal{Q}^-_2$ required a similar number of DOFs as the trimmed serendipity elements $\mathcal{S}^-_3$.  Compared this way, the trimmed serendipity elements provide an extra order of magnitude of accuracy over the tensor product element.  Furthermore, in \cref{fig:3dMixedTimeError} the time required for $\mathcal{S}^-_3$ and $\mathcal{Q}^-_2$ was also approximately equal.   Thus while it is helpful to compare $\mathcal{Q}^-_r$ and $\mathcal{S}^-_r$ to see that the trimmed serendipity elements have the expected convergence behavior, a more practical computational comparison is between $\mathcal{Q}^-_r$ and $\mathcal{S}_{r+1}^-$.  

The eigenvalue problem yields another example of comparing the tensor product and trimmed serendipity elements, where instead of refining the mesh, we  refined the order of the element used.  Just as in the mixed Poisson problem, we again see that \cref{fig:EigenvalueDofsError} shows $\mathcal{S}^-_2$ has a higher error for $\lambda = 3$ than $\mathcal{Q}^-_2$.  However comparing against where the DOFs are approximately equal leads to a comparison between $\mathcal{S}^-_3$ and $\mathcal{Q}^-_2$.  In this scenario, we had that $\mathcal{Q}^-_2$ required 104544 DOFs yielding an error of $1.33 \times 10^{-4}$ for $\lambda = 5$ while $\mathcal{S}^-_3$ required 106896 and achieved an error of $3.67 \times 10^{-7}$ for $\lambda = 5$.  In this case, we note that the time required for $\mathcal{S}^-_3$ did require more time to solve, using about 2.71 seconds while the $\mathcal{Q}^-_2$ required 1.98 seconds.

Our computational findings suggest that trimmed serendipity elements could be particularly beneficial at improving accuracy for compute-bound applications.  For any application, there is eventually a mesh resolution and element order for which refining the mesh or increasing the tensor product order is computationally infeasible.  In this instance, keeping the mesh but switching to a trimmed serendipity method of one order higher presents a new option to the practitioner that still provides an increase in accuracy without a significant increase to computational cost.

\newpage
\section*{Code availability}

All major Firedrake components, as well as the code for the numerical experiments in the paper have been archived on \citet{zenodo/Firedrake-20210419.2}. 

\begin{acks}
This work was supported by National Science Foundation (NSF) Collaborative Research Awards DMS-1913094 and DMS-1912653 and by U.S. Department of Energy, Office of Science, Office of Advanced Scientific Computing Research, under Award Number DE-SC-0019039. 
\end{acks}

\section{Appendix: Solver Configurations}
The solver configurations for the primal and mixed Poisson formulations can be found below.
\begin{lstlisting}[float=htbp,caption={An example of some solver parameters that we can use for the Poisson problem.  The options presented here solve the algebraic system with a simplified Newton method where the Jacobian is held constant at the first iterate.  Therefore it is factored at the beginning and triangular solves are applied to it at each subsequent iteration.  This has the effect of performing iterative refinement \cite{wilkinson1994rounding,moler1967iterative} and yields an increased algebraic accuracy on fine meshes.}, label={lst:solver_parameters}, numbers=left, firstnumber=1, xleftmargin=20pt,  xrightmargin=20pt]
...
params = {"snes_type": "newtonls",
          "snes_linesearch_type": "basic",
          "snes_monitor": None,
          "snes_converged_reason": None,
          "mat_type": "aij",
          "snes_max_it": 10,
          "snes_lag_jacobian": -2,
          "snes_lag_preconditioner": -2,
          "ksp_type": "preonly",
          "ksp_converged_reason": None,
          "ksp_monitor_true_residual": None,
          "pc_type": "lu",
          "snes_rtol": 1e-12,
          "snes_atol": 1e-20,
          "pc_factor_mat_solver_type": "mumps",
          "mat_mumps_icntl_14": "200",
          "mat_mumps_icntl_11": "2"}
...
\end{lstlisting}

\bibliographystyle{ACM-Reference-Format}
\bibliography{serendipityFiredrakePaper}
\end{document}